\documentclass[12pt]{article}

\input{preamble.tex}

\title{\large\bfseries\MakeUppercase{Quantum modularity and resurgence of the} \texorpdfstring{$q$}{q}-\MakeUppercase{Pochhammer symbol}}

\author{
	Murtadha Aljanabi%
	\texorpdfstring{\\ \normalsize\href{mailto:maljan1@lsu.edu}{maljan1@lsu.edu}}{}
	\and
	Koustav Mondal%
	\texorpdfstring{\footnote{\url{https://sites.google.com/view/kmondal}} \\ \normalsize \href{mailto:kmondal@southalabama.edu}{kmondal@southalabama.edu}}{}
}

\date{November 25, 2026}



\begin{document}
\maketitle

\begin{abstract}
	\noindent We investigate the resurgent and quantum modular properties of the logarithm of the $q$-Pochhammer symbol. From its asymptotic expansion, we compute the associated Stokes data and derive explicit formulas for its Borel--Laplace and median resummations. We show that the Borel--Laplace sums naturally encode the obstruction to modularity, leading to a quantum modular interpretation compatible with the modular resurgence paradigm of Fantini and Rella. Our results unify the two families of $q$-Pochhammer symbols studied by Fantini and Rella within a single framework, extending their results and providing further evidence for the interplay between resurgence and quantum modularity.
\end{abstract}

\pdfbookmark[1]{Contents}{toc}
\tableofcontents

\section{Introduction}
The $q$-Pochhammer symbol occupies a fundamental position in number theory, combinatorics, the theory of special functions, and mathematical physics. It plays a crucial role in the theory of basic hypergeometric series and partition identities, and it appears naturally in the study of modular forms, quantum invariants, and statistical mechanics \cite{andrews1974applications, andrews1986q, gasper2004basic}. In particular, the logarithm of a $q$-Pochhammer symbol provides a natural setting in which the interplay between modularity, asymptotic analysis, and resurgence is a rich topic of interest.

In this paper we study the logarithm of the $q$-Pochhammer symbol, which we define for $\tau$ in the complex upper-half plane $\H$ and $\r = \prmatrix*{r_1}{r_2} \in \R^2 \sm (\Z \times \Z_{\le 0})$ as
\begin{equation} \label{eq:fr-def}
	f_\r (\tau)
	= \slog\qpoch{\e(r_2 \tau - r_1)}{\e(\tau)}
	\coloneq \sum_{n=0}^{\infty} \log(1 - \e((n + r_2) \tau - r_1)),
\end{equation}
where $\qpoch{a}{q}$ denotes the $q$-Pochhammer symbol, $\e(z) \coloneqq e^{2\pi iz}$, and $\slog$ denotes the branch of the logarithm determined by the sum on the right-hand side with $\log$ being the principal branch. Our principal objective is to study the asymptotic expansion of $f_\r$, which we denote by $\tf_\r$, and the extent to which its resummation recovers the original function and its modular transformation properties.

The study of functions exhibiting modular behavior at rational points has a long history. An important precursor to the theory of quantum modular forms is the work of Lawrence and Zagier \cite{lawrence1999modular}, which connected modular forms and quantum invariants of three-manifolds through functions whose boundary values exhibit a subtle form of modularity. Zagier subsequently introduced the term \emph{quantum modular forms} in \cite{zagier2010quantum} to describe functions defined on rational points whose failure to transform modularly is nevertheless controlled by a suitably regular function. Unlike classical modular forms, quantum modular forms need not be holomorphic functions on the upper half-plane; instead, their modular properties are encoded in the regularity of their transformation obstructions. This perspective has since led to connections with Eichler integrals, mock modular forms, partial theta functions, knot invariants, and other arithmetic and geometric phenomena \cite{bringmann2020quantum, bringmann2016half, kopp2024shintani, mcspirit2025quantum}.

A complementary perspective is provided by \'Ecalle's theory of resurgence, developed to analyze divergent asymptotic expansions and the analytic information encoded in them. In resurgence theory, the Borel transform of a divergent formal series converts its coefficients into an analytic function whose singularities govern the exponentially small contributions that are invisible to ordinary asymptotic expansions. Borel--Laplace summation reconstructs analytic functions from suitable asymptotic data, while the Stokes phenomenon describes the changes in these sums as the summation direction crosses singular directions in the Borel plane. Median resummation, in particular, provides a prescription for treating such directions when lateral Borel sums differ. These ideas have found applications across differential equations, mathematical physics, and the study of perturbative and non-perturbative expansions.

The relationship between resurgence and quantum modularity has recently attracted substantial attention. In \cite{fantini2024modular}, Fantini and Rella proposed the framework of \emph{modular resurgence}, conjecturing that suitable arithmetic structures in the Borel plane, particularly those reflected in the Stokes constants and their generating functions, are related to quantum modularity. Their conjectures predict, under appropriate hypotheses, that median resummation reconstructs the original $q$-series and that the resulting functions exhibit quantum modular behavior. In a related development, Fantini and Rella established resurgent and quantum modular properties for functions arising from topological string theory in \cite{fantini2024modular}. McSpirit and Rolen \cite{mcspirit2025quantum} subsequently investigated the relationship between quantum modular forms and resurgence, demonstrating how asymptotic expansions and their resummations can unify and strengthen several examples of quantum modularity. More recently, Fantini and Rella \cite{fantini2025modular} studied the summability and quantum modularity of character-weighted sums of $q$-Pochhammer symbols. Together, these works provide a growing body of evidence that resurgence can furnish an analytic framework for understanding quantum modular phenomena.

Before stating our results, we introduce the following notation to be used throughout the paper.

\begin{definition}
	For $\r = \prmatrix*{r_1}{r_2} \in \R^2$, let $\{\r\} \coloneqq \prmatrix*{r_1 - \floor{r_1}}{r_2 - \ceil{r_2} + 1}$.
\end{definition}

Note that $\{\r\} \equiv \r \pmod{\Z^2}$ and $\{\r\} \in [0,1) \times (0,1]$. If $\{x\} \coloneqq x - \floor{x}$ denotes the fractional part of a real number $x$, then $\{\r\} = \prmatrix*{\{r_1\}}{1 - \{-r_2\}}$. Equivalently, $\{\r\} = \prmatrix*{\{r_1\}}{\{r_2\}}$ unless $r_2 \in \Z$, in which case $\{\r\} = \prmatrix*{\{r_1\}}{1}$.

\paragraph{Main results}
The purpose of this paper is to investigate the preceding ideas for the family of logarithmic $q$-Pochhammer symbols $f_\r$. Using the asymptotic expansion $\tf_\r$, we study the associated Borel transforms and Stokes phenomena to unveil the resurgent structure and relate that to the obstruction to modularity for $f_\r$. We also use these results to analyze median resummation of $\tf_\r$. We denote by $\S_\theta (\tf_\r)$ the Borel--Laplace sum of $\tf_\r$ in the direction $\theta$, and by $\Smed_\theta (\tf_\r)$ its median resummation. The precise definitions of these operators are given in \zcref{sec:resurge}. Explicitly, we have the following results.

\begin{theorem}[\zcref{thm:bl-fr-r1-nonint,thm:bl-fr-r1-int}] \label{thm:bl-fr}
	Let $\r = \prmatrix*{r_1}{r_2} \in (\R \times \R_{\ge 0}) \sm (\Z \times \{0\})$ and $\tau \in \H$. The Borel--Laplace sums of the asymptotic series of $f_\r (\tau)$ as $\tau \to 0$ along the rays $\theta = 0$ and $\theta = \pi$ satisfy
	\begin{alignat}{2}
		\S_0   (\tf_\r) (\tau) &= f_\r (\tau) - f_{\{ S\r\}} (S \cdot \tau), & \qquad &\Re(\tau) > 0, \\
		\S_\pi (\tf_\r) (\tau) &= f_\r (\tau) - f_{\{-S\r\}} (S \cdot \tau), & \qquad &\Re(\tau) < 0,
	\end{alignat}
	where $S \coloneqq \prmatrix*{0 & -1}{1 & 0}$ and $S \cdot \tau \coloneqq -1/\tau$.
\end{theorem}

In particular, \zcref{thm:bl-fr} proves Conjecture~5 of Fantini--Rella \cite{fantini2025modular}, which corresponds to the special case $\r = \prmatrix*{0}{k/N}$. We give a precise identification in \zcref{rmk:FR-conj}. The following result is a direct consequence of \zcref{thm:bl-fr}.

\newcommand{\frmedresum}{
	Let $\r = \prmatrix*{r_1}{r_2} \in (\R \times \R_{\ge 0}) \sm (\Z \times \{0\})$ and $\tau \in \H$. Then the median resummation of $\tf_\r$ along the ray $\theta = \pi/2$ is given by
	\begin{equation}
		\Smed_{\pi/2} (\tf_\r) (\tau)
		= f_\r (\tau) - \frac{1}{2} \left(f_{\{S\r\}} (S \cdot \tau) + f_{\{-S\r\}} (S \cdot \tau)\right).
	\end{equation}
}

\begin{corollary}[\zcref{cor:fr-med-resum}] \label{cor:fr-med-resum-intro}
	\frmedresum
\end{corollary}

As seen in \zcref{cor:fr-med-resum-intro}, the median resummation of $\tf_\r$ generally differs from $f_\r$ by a correction term. Replacing $\r$ by $\{-\r\}$ interchanges the two summands in this correction term, so it cancels in the difference $f_\r - f_{\{-\r\}}$. The median resummation therefore recovers this difference exactly.

\newcommand{\Frmedresum}{
	For $\r = \prmatrix*{r_1}{r_2} \in (\R \times \R_{\ge 0}) \sm (\Z \times \{0\})$ and $\tau \in \H$, let
	\begin{alignat}{2}
		F_\r (\tau)                   &= f_\r (\tau) - f_{\{-\r\}} (\tau).          &        &           \\ \shortintertext{Then}
		\S_0 (\tF_\r) (\tau)          &= F_\r (\tau) - F_{\{S\r\}} (S \cdot \tau),  & \qquad & \Re(\tau) > 0,                  \\
		\S_\pi (\tF_\r) (\tau)        &= F_\r (\tau) - F_{\{-S\r\}} (S \cdot \tau), & \qquad & \Re(\tau) < 0, \quad \text{and} \\
		\Smed_{\pi/2} (\tF_\r) (\tau) &= F_\r (\tau).                               &        &
	\end{alignat}
}

\begin{corollary}[\zcref{cor:Fr-med-resum}]
	\Frmedresum
\end{corollary}

\paragraph{Outline of the paper}
The paper is organized as follows. In \zcref{sec:background}, we review the necessary background on the $q$-Pochhammer symbols, including its asymptotic expansion, and recall the relevant notions of quantum modularity and modular resurgence. In \zcref{sec:fr-resurgence}, we develop the resurgent analysis of $f_\r$ using the asymptotic expansions, and calculate Borel singularities, Stokes phenomena, and Borel--Laplace sums. Finally in \zcref{sec:fr-quantum} we establish the quantum modularity of $f_\r$ through the Borel--Laplace sum of $\tf_\r$.

\paragraph{Acknowledgments}
The authors express their gratitude to Gene S.~Kopp for his invaluable guidance and insightful discussions. They also gratefully acknowledge Veronica Fantini, Amanda Folsom, Eleanor McSpirit, Claudia Rella, and Campbell Wheeler for helpful discussions.

The first author was partially supported by NSF grants DMS-2302514 and DMS-2601235. The second author was partially supported by NSF grant DMS-2302514.

\section{Background} \label{sec:background}

\subsection{The \texorpdfstring{$q$}{q}-Pochhammer symbol} \label{sec:qpoch}
For $a, q \in \C$ and $n \in \Z_{>0}$, the \emph{$q$-Pochhammer symbol} is defined as
\begin{equation}
	\qpoch{a}{q}[n] = \prod_{k=0}^{n - 1} (1 - aq^k).
\end{equation}
If, in addition, we have that $\abs{q} < 1$, then the product converges absolutely and locally uniformly as $n \to \infty$, so we set
\begin{equation}
	\qpoch{a}{q} \coloneqq \lim_{n \to \infty} \qpoch{a}{q}[n] = \prod_{k=0}^{\infty} (1 - aq^k).
\end{equation}
Further, if $\abs{a} \le 1$ and $a \ne 1$, then $\log(1 - aq^k)$ is continuous in $q$ for every $k \ge 0$ and the series
\begin{equation}
	\sum_{k=0}^{\infty} \log(1 - aq^k)
\end{equation}
converges uniformly on compact subsets of the unit disk $\{q \in \C \colon \abs{q} < 1\}$. It follows that the function
\begin{equation}
	\slog\qpoch{a}{q} \coloneqq \sum_{k=0}^{\infty} \log(1 - aq^k)
\end{equation}
is a continuous choice of logarithm of $\qpoch{a}{q}$ for $\abs{q} < 1$.

Now set $q = \e(\tau)$ and $a = \e(r_2 \tau - r_1)$. Then for $\r = \prmatrix*{r_1}{r_2} \in (\R \times \R_{\ge 0}) \sm (\Z \times \{0\})$, and $\tau \in \H$, we have $\abs{q} < 1$ and $\abs{a} \le 1$ with $a \ne 1$. Hence the function
\begin{equation}
	f_\r (\tau) = \slog\qpoch{\e(r_2 \tau - r_1)}{\e(\tau)}
\end{equation}
is continuous (and holomorphic) in $\tau$.

Our results rely on an asymptotic expansion of $f_\r (\tau)$ near $\tau = 0$ and an integral representation of a certain ratio of $q$-Pochhammer symbols. We state these required results below.

\begin{lemma}[{\cite[p.~178]{wheeler2022modular}, \cite[Thm.~1]{katsurada2003asymptotic}}] \label{lem:fr-asymp}
	Let $\r = \prmatrix*{r_1}{r_2} \in (\R \times \R_{\ge 0}) \sm (\Z \times \{0\})$ and let $\tf_\r (\tau)$ denote the asymptotic of $f_\r (\tau)$ as $\tau \to 0$. If $r_1 \notin \Z$, then as $\tau \to 0$ with $\delta < \arg(\tau) < \frac{\pi}{2} - \delta$ for some $\delta > 0$,
	\begin{equation}
		\tf_\r (\tau)
		= \frac{1}{2\pi i \tau} \Li_2 (\e(-r_1)) - B_1 (r_2) \log(1 - \e(-r_1)) + \psi_\r (\tau),
	\end{equation}
	where
	\begin{equation}
		\psi_\r (\tau) = \sum_{k=0}^{\infty} \frac{B_{k + 2} (r_2) \Li_{-k} (\e(-r_1))}{(k + 2)!}\, (2\pi i \tau)^{k + 1}.
	\end{equation}
	On the other hand, if $r_1 \in \Z$, then as $\tau \to 0$ with $\delta < \arg(\tau) < \pi - \delta$,
	\begin{equation}
		\tf_\r (\tau)
		= -\frac{\pi i}{12\tau} - B_1 (r_2) \log(-2\pi i \tau) - \log\frac{\Gamma(r_2)}{\sqrt{2\pi}} + \psi_{r_2} (\tau),
	\end{equation}
	where
	\begin{equation}
		\psi_{r_2} (\tau)
		= -\sum_{k=1}^{\infty} \frac{B_k (1) B_{k + 1} (r_2)}{k (k + 1)!}\, (2\pi i \tau)^k.
	\end{equation}
\end{lemma}

\begin{remark} \label{rmk:fr-asymp}
	The expansion $\tf_\r (\tau)$ for $r_1 \notin \Z$ is also valid as $\tau \to 0$ with $\frac{\pi}{2} + \delta < \arg(\tau) < \pi - \delta$ because $f_\r (\tau) = \ol{f_{\prmatrix*{-r_1}{r_2}} (-\ol{\tau})}$.
\end{remark}

\begin{theorem}[{\cite[p.~182]{wheeler2022modular}}] \label{thm:wheelerqpochhammer}
	Let $\tau \in \H$ with $\Re(\tau) > 0$. If $z \in \C \sm \Z$ and $M \in \Z$ such that $M\tau + z \notin \Z$, then
	\begin{equation}
		\slog\frac{\qpoch{\e(z)}{\e(\tau)}[M + 1]}
		{\qpoch{\e((z - \ceil{\Re(z)}) / \tau)}{\e(-1/\tau)}[N]}
		= \varphi(\tau, z, M),
	\end{equation}
	where the logarithm on the left-hand side is defined as
	\begin{equation}
		\sum_{n=0}^{M} \log(1 - \e(n \tau + z))
		- \sum_{n=0}^{N - 1} \log(1 - \e((z - \ceil{\Re(z)} - n) / \tau))
	\end{equation}
	with $N = \ceil{\Re(M\tau + z)} - \ceil{\Re(z)}$, and
	\begin{align}
		\varphi(\tau, z, M)
		= \frac{1}{2\pi i\tau}& \Li _2(\e(z)) - \frac{1}{2\pi i\tau} \Li_2 (\e(M\tau + z)) \\
		&+ \frac{1}{2} \log(1 - \e(z)) + \frac{1}{2} \log(1 - \e(M\tau + z)) \\
		&+ \frac{i}{\tau} \int_{0}^{\infty} \frac{\log(1 - \e(it + z))
		- \log(1 - \e(-it + z))}{\e(-it/\tau) - 1}\, dt \\
		&- \frac{i}{\tau} \int_{0}^{\infty} \frac{\log(1 - \e(M\tau + it + z))
		- \log(1 - \e(M\tau - it + z))}{\e(-it/\tau) - 1}\, dt.
	\end{align}
\end{theorem}

The leading logarithmic terms appearing in $\varphi(\tau, z, M)$ can be absorbed into the integrals by a suitable change of variables followed by integration by parts, which we carry out in \zcref{cor:int-to-qpoch} after sending $M \to \infty$. Before stating and proving the \namecref{cor:int-to-qpoch}, we establish the following result, which fixes a continuous choice of branch of $\log(1 - e^z)$ that will be used throughout the paper.

\begin{proposition} \label{prop:continuous-L}
	For $z \in \C$, define
	\begin{equation}
		\ell(z) =
		\begin{cases}
			\log(1 - e^z),                & \Re(z) \le 0, \\
			z - \pi i + \log(1 - e^{-z}), & \Re(z) > 0.
		\end{cases}
	\end{equation}
	Then $\ell$ is a continuous choice of branch of $\log(1 - e^z)$ on $\C \sm \{iy \colon y \le 0 \text{ or } y \ge 2\pi\}$.
\end{proposition}

\begin{proof}
	Observe that for $\Re(z) > 0$,
	\begin{equation}
		e^{\ell(z)}
		= e^{z - \pi i + \log(1 - e^{-z})}
		= -e^{z} (1 - e^{-z})
		= 1 - e^z
		= e^{\log(1 - e^z)},
	\end{equation}
	so $\ell(z)$ is a choice of branch of $\log(1 - e^z)$ for $\Re(z) > 0$.
	
	If $\Re(z) < 0$, then $\abs{e^z} < 1$, so $\ell(z) = \log(1 - e^z)$ is continuous. Similarly, if $\Re(z) > 0$, then $\abs{e^{-z}} < 1$ and $\log(1 - e^{-z})$ is continuous, which implies that $\ell(z)$ is continuous. It remains to show that $\ell(z)$ is continuous at $z = iy$ for $y \in (0, 2\pi)$. Indeed, since
	\begin{equation}
		1 - e^{iy}
		= -e^{iy/2} (e^{iy/2} - e^{-iy/2})
		= -2i e^{iy/2} \sin(y/2)
		= 2 \sin(y/2) e^{i(y - \pi)/2},
	\end{equation}
	it follows that
	\begin{equation}
		\ell(iy)
		= \log(1 - e^{iy})
		= \log(2 \sin(y/2)) + i\, \frac{y - \pi}{2}.
	\end{equation}
	On the other hand, since $1 - e^{-iy} = \ol{1 - e^{iy}}$, we have
	\begin{align}
		\lim_{\varepsilon \to 0^+} \ell(\varepsilon + iy)
		&= iy - \pi i + \log(1 - e^{-iy}) \\
		&= i(y - \pi) + \log(2 \sin(y/2)) - i\, \frac{y - \pi}{2} \\
		&= \log(2 \sin(y/2)) + i\, \frac{y - \pi}{2} \\
		&= \ell(iy).
	\end{align}
	This proves that $\ell(z)$ is continuous at $z = iy$, which completes the proof.
\end{proof}

\begin{corollary} \label{cor:int-to-qpoch}
	Let $\tau \in \H$ with $\Re(\tau) > 0$ and let $z \in \C$ such that $\Re(z) \notin \Z$. Then
	\begin{equation}
		-\frac{1}{2\pi i} \int_{-\infty}^{\infty} \frac{\ell(t/\tau)}{1 - \e(-z) e^{t}}\, dt
		= \slog\frac{\qpoch{\e(z)}{\e(\tau)}}{\qpoch{\e((z - \ceil{\Re(z)}) / \tau)}{\e(-1/\tau)}},
	\end{equation}
	where $\ell(z)$ is the logarithm of $1 - e^z$ defined in \zcref{prop:continuous-L}.
\end{corollary}

\begin{proof}
	By sending $M \to \infty$ in \zcref{thm:wheelerqpochhammer}, we get
	\begin{equation} \label{eq:log-qpoch-lim-psi}
		\slog\frac{\qpoch{\e(z)}{\e(\tau)}}{\qpoch{\e((z - \ceil{\Re(z)}) / \tau)}{\e(-1/\tau)}}
		= \lim_{M \to \infty} \varphi(\tau, z, M),
	\end{equation}
	with
	\begin{equation} \label{eq:lim-psi}
		\lim_{M \to \infty} \varphi(\tau, z, M)
		= \frac{1}{2\pi i\tau} \Li _2(\e(z)) + \frac{1}{2} \log(1 - \e(z))
		+ \frac{i}{\tau} \int_{0}^{\infty} \frac{g(t) - g(-t)}{\e(-it/\tau) - 1}\, dt,
	\end{equation}
	where $g(t) = \log(1 - \e(it + z))$. For $\varepsilon > 0$, let
	\begin{equation}
		I_\varepsilon
		= \int_{\varepsilon}^{\infty} \frac{g(t) - g(-t)}{\e(-it/\tau) - 1}\, dt
		= \int_{-\infty}^{-\varepsilon} \frac{g(t)}{1 - \e(it/\tau)}\, dt + \int_{\varepsilon}^{\infty} \frac{g(t)}{\e(-it/\tau) - 1}\, dt.
	\end{equation}
	Using $1/(e^{-\alpha} - 1) = 1/(1 - e^\alpha) - 1$ in the second integral on the right-hand side and sending $\varepsilon \to 0^+$ yields
	\begin{align}
		\lim_{\varepsilon \to 0^+} I_\varepsilon
		&= \lim_{\varepsilon \to 0^+} \left(\int_{-\infty}^{-\varepsilon} \frac{g(t)}{1 - \e(it/\tau)}\, dt
		+ \int_{\varepsilon}^{\infty} \frac{g(t)}{1 - \e(it/\tau)}\, dt - \int_{\varepsilon}^{\infty} g(t)\, dt\right) \\
		&= \PV \int_{-\infty}^{\infty} \frac{g(t)}{1 - \e(it/\tau)}\, dt - \int_{0}^{\infty} g(t)\, dt,
	\end{align}
	where $\PV$ denotes the Cauchy principal value (which is necessary due to the pole at $t = 0$). Using the substitution $u = \e(it + z)$, we find that
	\begin{equation}
		\int_{0}^{\infty} g(t)\, dt
		= \frac{1}{2\pi} \int_{0}^{\e(z)} \frac{\log(1 - u)}{u}\, du
		= -\frac{1}{2\pi} \Li_2 (\e(z)).
	\end{equation}
	Substituting into \eqref{eq:lim-psi} then replacing $t$ by $t/(2\pi)$ in the integral yields
	\begin{align}
		\lim_{M \to \infty} \varphi(\tau, z, M)
		&= \frac{1}{2} \log(1 - \e(z)) + \frac{i}{\tau} \PV \int_{-\infty}^{\infty} \frac{g(t)}{1 - \e(it/\tau)}\, dt \\
		&= \frac{1}{2} \log(1 - \e(z)) - \frac{1}{2\pi i \tau} \PV \int_{-\infty}^{\infty} \frac{g(t/(2\pi))}{1 - e^{-t/\tau}}\, dt
		\intertext{which, by integrating by parts with $\tau \cdot d(\ell(t/\tau)) = 1/(1 - e^{-t/\tau})\, dt$, becomes}
		&= \frac{1}{2} \log(1 - \e(z)) - \frac{1}{2\pi i} \PV \int_{-\infty}^{\infty} g(t/(2\pi))\, d(\ell(t/\tau)) \\
		&= \frac{1}{2} \log(1 - \e(z)) - \frac{1}{2\pi i} B - \frac{1}{2\pi i} \int_{-\infty}^{\infty} \frac{\ell(t/\tau)}{1 - \e(-z) e^t}\, dt,
	\end{align}
	where $B$ is the sum of the boundary terms, given by
	\begin{align}
		B
		&= \lim_{\varepsilon \to 0^+} \left(g(t/(2\pi)) \ell(t/\tau) \big|_{t=-\infty}^{t=-\varepsilon}
		+ g(t/(2\pi)) \ell(t/\tau) \big|_{t=\varepsilon}^{t=\infty}\right) \\
		&= g(0) \lim_{\varepsilon \to 0^+} \bigl(\ell(-\varepsilon/\tau) - \ell(\varepsilon/\tau)\bigr) \\
		&= g(0) \lim_{\varepsilon \to 0^+} \bigl(\log(1 - e^{-\varepsilon/\tau}) - \varepsilon/\tau + \pi i - \log(1 - e^{-\varepsilon/\tau})\bigr) \\
		&= \pi i \log(1 - \e(z)).
	\end{align}
	Hence,
	\begin{equation}
		\lim_{M \to \infty} \varphi(\tau, z, M)
		= -\frac{1}{2\pi i} \int_{-\infty}^{\infty} \frac{\ell(t/\tau)}{1 - \e(-z) e^t}\, dt.
	\end{equation}
	Substituting into \eqref{eq:log-qpoch-lim-psi} yields the desired result.
\end{proof}

\begin{remark} \label{rmk:int-to-qpoch}
	For $\Re(\tau) < 0$, we use \zcref{cor:int-to-qpoch} with $(z, \tau)$ replaced by $(-\bar{z}, -\bar{\tau})$ then take the complex conjugate of both sides to get
	\begin{equation}
		\frac{1}{2\pi i} \int_{-\infty}^{\infty} \frac{\ell(-t/\tau)}{1 - \e(-z) e^{t}}\, dt + \log(1 - \e(z))
		= \slog\frac{\qpoch{\e(z)}{\e(\tau)}}{\qpoch{\e((\floor{\Re(z)} - z) / \tau)}{\e(-1/\tau)}},
	\end{equation}
	with the additional logarithmic term on the left-hand side appearing because $\Re(-t/\bar{\tau}) > 0$ for $t > 0$, which implies that
	\begin{equation}
		\ol{\ell(-t/\bar{\tau})}
		= \ol{-t/\bar{\tau} - \pi i + \log(1 - e^{t/\bar{\tau}})}
		= -t/\tau + \pi i + \log(1 - e^{t/\tau})
		= \ell(-t/\tau) + 2\pi i.
	\end{equation}
\end{remark}

\begin{remark}
	The condition $\Re(z) \notin \Z$ in \zcref{cor:int-to-qpoch} is necessary because the integrand has a pole when $\Re(z) \in \Z$ and so the integral diverges. In fact, if $\Re(z) \in \Z$, then for $\Re(\tau) > 0$,
	\begin{equation}
		-\frac{1}{2\pi i} \PV \int_{-\infty}^{\infty} \frac{\ell(t/\tau)}{1 - \e(-z) e^{t}}\, dt
		- \frac{1}{2} \log(1 - \e(i \Im(z) / \tau))
		= \slog\frac{\qpoch{\e(z)}{\e(\tau)}}{\qpoch{\e(i \Im(z) / \tau)}{\e(-1/\tau)}},
	\end{equation}
	and for $\Re(\tau) < 0$,
	\begin{align}
		\frac{1}{2\pi i} \PV \int_{-\infty}^{\infty} \frac{\ell(-t/\tau)}{1 - \e(-z) e^{t}}\, dt + \log(1 - \e(z))
		- \frac{1}{2} \log(1 - \e(-i \Im(z) / \tau)) \\
		= \slog\frac{\qpoch{\e(z)}{\e(\tau)}}{\qpoch{\e(-i \Im(z) / \tau)}{\e(-1/\tau)}}.
	\end{align}
\end{remark}

\subsection{Faddeev's quantum dilogarithm}
In this section, we briefly review Faddeev's quantum dilogarithm and its properties, including its meromorphic continuation and expression in terms of $q$-Pochhammer symbols \cite{faddeev1994quantum, garoufalidis2021resurgence}. We begin by formally defining it.

\begin{definition}[Faddeev \cite{faddeev1995discrete}] \label{def:faddeev}
	Let
	\begin{equation}
		c_b = \frac{i}{2}(b + b^{-1}).
	\end{equation}
	For $\abs{\Im(z)} < \abs{\Im(c_b)}$, define \emph{Faddeev's quantum dilogarithm} $\Phi_b (z)$ by the integral
	\begin{equation} \label{eq:faddeevintegral}
		\Phi_b (z) = \exp\left(\int_{\R + i\varepsilon} \frac{e^{-2 itz}}{4\sinh(tb) \sinh (t b^{-1})}\, \frac{dt}{t}\right).
	\end{equation}
\end{definition}

\begin{remark} \label{rmk:faddeevcont}
	It is clear from the definition that $\Phi_b (z)$ satisfies the symmetric properties
	\begin{equation}
		\Phi_b (z) = \Phi_{-b} (z) = \Phi_{b^{-1}} (z).
	\end{equation}
	By means of this symmetry, \eqref{eq:faddeevintegral} can be extended to the region $\Im(b^2) < 0$ and it also admits an analytic continuation to all values of $b$ where $b^2 \notin \R_{\le 0}$.
	
	Moreover, as a function of $z$, $\Phi_b (z)$ has a meromorphic continuation to the complex plane with
	\begin{enumerate}
		\item an essential singularity at infinity,
		\item poles at $z = c_b + imb + in b^{-1}$ for $m,n \in \Z_{>0}$,
		\item zeroes at $z = -c_b - imb - in b^{-1}$ for $m,n \in \Z_{>0}$.
	\end{enumerate}
\end{remark}

\begin{lemma} \label{lem:faddeev-qpoch}
	When $\Im(b^2) > 0$, the integral formula of $\Phi_b (z)$ in \eqref{eq:faddeevintegral} can be expressed in terms of a ratio of $q$-Pochhammer symbols as
	\begin{equation}
		\Phi_b (z) = \frac{\qpoch{e^{2\pi b(z + c_b)}}{\e(b^2)}}{\qpoch{e^{2\pi b^{-1}(z - c_b)}}{\e(-b^{-2})}}.
	\end{equation}
\end{lemma}

\begin{proof}
	See \cite[Prop.~A.3]{andersen2016asymptotic}.
\end{proof}

By substituting $z = i r_2 b - i(r_1 - 1) b^{-1} - c_b$ in \zcref{lem:faddeev-qpoch} with $\tau = b^2 \in \H$ and $\r = \prmatrix*{r_1}{r_2} \in \R_{>0}^2$, we see that Faddeev's quantum dilogarithm relates to $f_\r (\tau)$ with $b^2 = \tau$ as
\begin{align}
	\slog {} \Phi_b (z)
	&=
	\begin{multlined}[t][0.8\textwidth]
		\slog\qpoch[\big]{\exp(2\pi i r_2 \tau - 2\pi i (r_1 - 1)}{\e(\tau)} \\
		- \slog\qpoch[\big]{\exp(2\pi i r_2 - 2\pi i (r_1 - 1) b^{-2} - 4\pi b^{-1} c_b)}{\e(-1/\tau)}
	\end{multlined} \\
	&= \slog\qpoch{\e(r_2 \tau - r_1)}{\e(\tau)} - \slog\qpoch{\e(r_2 - r_1 / \tau)}{\e(-1/\tau)} \\
	&= f_\r (\tau) - f_{S\r} (S \cdot \tau).
\end{align}

By \zcref{rmk:faddeevcont}, the zeros and poles of $\Phi_b (z)$ are at $z = \pm (c_b + imb + in b^{-1})$ for $m,n \in \Z_{>0}$. Therefore, the singularities of $f_\r (\tau) - f_{S\r} (S \cdot \tau)$ occur when
\begin{equation}
	z
	= i r_2 b - i(r_1 - 1) b^{-1} - c_b
	= \pm (c_b + imb + in b^{-1}).
\end{equation}
Multiplying by $b$ then using $b^2 = \tau$ and solving for $\tau$ yields
\begin{equation}
	\tau = \frac{r_1 + n}{r_2 - m - 1}
	\qqor
	\tau = \frac{r_1 - n - 1}{r_2 + m}.
\end{equation}

\subsection{Quantum modular forms}
The notions of ``quantum modular forms'' and ``holomorphic quantum modular forms'' are purposely intended (by Zagier) to be somewhat vague, so as to be expandable to fit many interesting examples \cite{zagier2010quantum}. For our purposes, we use a more precise definition, also due to Zagier \cite{zagier2020holomorphic}. If $f$ is a function on some complex domain $U \subseteq \C$ and $k \in \Z$, the \emph{slash operator of weight $k$} defines for each $\gamma = \prmatrix*{a & b}{c & d} \in \SL_2 (\Z)$ another function
\begin{equation}
	(f|_k\gamma)(\tau) \coloneqq j_\gamma(\tau)^{-k} f(\gamma\cdot\tau),
\end{equation}
where $j_\gamma(\tau) \coloneqq c \tau + d$ and $\gamma \cdot \tau \coloneqq \frac{a\tau + b}{c\tau + d}$. The function $f|_k \gamma$ is defined on $\gamma^{-1} \cdot U \sm \{-d/c\}$ if $k > 0$ and on $\gamma^{-1} \cdot U$ if $k \le 0$. If $k = 0$, we drop the subscript and write $f | \gamma$ instead of $f |_k \gamma$.

\begin{definition} \label{def:holomorphic-qmf}
	A holomorphic function $f \colon \H \to \C$ is said to be a \emph{holomorphic quantum modular form} of weight $k \in \Z$ for a subgroup $\Gamma$ of $\SL_2 (\Z)$ if for all $\gamma = \prmatrix*{a & b}{c & d} \in \Gamma$, the obstruction to modularity
	\begin{equation}
		h_\gamma (\tau) \coloneqq f(\tau) - (f|_k \gamma)(\tau)
	\end{equation}
	extends to an analytic function on the cut plane $\C_\gamma \coloneqq \C \sm \{x \in \R \colon cx + d \le 0\}$.
\end{definition}

Note that it suffices to prove that $h_\gamma$ is analytic on $\C_\gamma$ for a set of generators of $\Gamma$, since
\begin{equation}
	h_{\gamma_1 \gamma_2} = h_{\gamma_1} |_k \gamma_2 + h_{\gamma_2}.
\end{equation}
If $\Gamma = \SL_2 (\Z)$, then it is known that $\Gamma = \langle S, T \rangle$ where $S = \prmatrix*{0 & -1}{1 & 0}$ and $T = \prmatrix*{1 & 1}{0 & 1}$.

\subsection{Resurgence} \label{sec:resurge}
Formal power series of the form $\psi(\tau) = \sum_{n=0}^{\infty} a_n \tau^n$ whose coefficients satisfy a factorial growth bound at large $n$ arise naturally in asymptotic analysis involving special functions, modular forms, and $q$-series. Although such expansions are generally divergent, their coefficients often encode rich analytic and arithmetic information. Resurgence theory provides a framework for extracting this information by studying the analytic properties of a transformed series. We begin by recalling the basic notions of Gevrey-$1$ series and Borel transformation before introducing the particular class of resurgent series relevant to this paper.

\paragraph{Gevrey-$1$ series and the Borel transform}
We say that the formal power series $\psi(\tau) = \sum_{n=0}^{\infty} a_n \tau^n$ is a \emph{Gevrey-$1$} series if there exists constants $M,k > 0$ such that
\begin{equation}
	\abs{a_n} < M k^n n!.
\end{equation}
Factorial growth of the coefficients means that the series is potentially divergent, but at a controlled rate. The Borel transform, defined below, compensates for this growth by dividing the coefficients by factorials, thereby producing a series with a positive radius of convergence. Since we are primarily interested in asymptotic expansions without a constant term, we use the following convention.

\begin{definition}[Borel transform]
	The \emph{Borel transform} is the linear map
	\begin{equation}
		\B \colon \tau\, \C[\mspace{-2mu}[\tau]\mspace{-2mu}] \to \C[\mspace{-2mu}[\xi]\mspace{-2mu}]
	\end{equation}
	defined on formal power series as
	\begin{equation}
		\B\left(\sum_{n=0}^{\infty} a_n \tau^{n + 1}\right) (\xi)
		= \sum_{n=0}^{\infty} \frac{a_n}{n!}\, \xi^n.
	\end{equation}
\end{definition}

If $\psi(\tau)$ is Gevrey-$1$, its Borel transform converges in a neighborhood of $\xi = 0$. The central question is whether this local analytic function can be continued beyond that neighborhood and, if so, what singularities obstruct its continuation. These singularities are fundamental to resurgence theory: they govern the additional analytic information encoded in a divergent asymptotic expansion.

\paragraph{Resurgent series}
A Gevrey-$1$ series $\psi(\tau)$ is called \emph{resurgent} if its Borel transform $\B(\psi)(\xi)$ admits analytic continuation along paths of arbitrary finite length, subject to avoiding an appropriate discrete set of singularities (see \cite[Def.~2.1]{fantini2025modular} for a precise definition). Moreover, we say that $\psi(\tau)$ is \emph{simple resurgent} if it is resurgent and $\B(\psi)(\xi)$ only has simple poles and logarithmic branch points.

For the remainder of this discussion, we restrict our attention to simple resurgent series whose Borel transforms have only simple poles. In this case, the analytic continuation of $\B(\psi)(\xi)$ is single-valued and defines a meromorphic function on the entire complex plane. Let $\{\xi_\omega\}_{\omega \in \Omega}$ denote its (countable) set of poles. The \emph{Stokes constant} associated with $\xi_\omega$ is defined by
\begin{equation}
	S_\omega \coloneqq -2\pi i \Res_{\xi = \xi_\omega} \B(\psi)(\xi).
\end{equation}

It is often the case that the poles of $\B(\psi) (\xi)$ lie on a single line through the origin, with the origin excluded. In this case, we re-index the poles by $m \in \Z_{\ne 0}$ so that those lying on one side of the line (with respect to the origin) are indexed by positive integers and those on the other side are indexed by negative integers. Thus, we can write
\begin{equation}
	\xi_m = A \rho_m, \qquad
	\xi_{-m} = -A \rho_{-m}, \qquad
	m \in \Z_{>0},
\end{equation}
where $A \in \C^\times$ and $\{\rho_{\pm m}\}_{m=1}^{\infty}$ are strictly increasing sequences of positive real numbers tending to infinity. We associate to the Stokes constants a pair of \emph{general Dirichlet series} (defined and discussed in \cite[Sec.~8.2]{apostol1990modular} and \cite[p.~65]{serre2012course}), which are defined as
\begin{equation} \label{eq:L-pm-general}
	L_+ (s) \coloneqq \sum_{m=1}^{\infty} \frac{S_m}{\rho_m^s}
	\qqand
	L_- (s) \coloneqq -\sum_{m=1}^{\infty} \frac{S_{-m}}{\rho_{-m}^s}.
\end{equation}

A particularly interesting situation occurs when the poles are evenly spaced along the line, that is,
\begin{equation}
	\xi_m = mA, \quad m \in \Z_{\ne 0}.
\end{equation}
In this case, $\rho_m = \rho_{-m} = m$, so that the general Dirichlet series $L_+ (s)$ and $L_- (s)$ become ordinary Dirichlet series:
\begin{equation} \label{eq:L-pm-ordinary}
	L_+(s) = \sum_{m=1}^{\infty} \frac{S_m}{m^s},
	\qquad
	L_-(s) = -\sum_{m=1}^{\infty} \frac{S_{-m}}{m^s}.
\end{equation}


\paragraph{Borel--Laplace summation}
The Borel transform converts a Gevrey-$1$ formal series into a convergent power series near the origin. When the resulting analytic function admits analytic continuation along a given direction with suitable growth, the Laplace transform provides an analytic function associated with the original formal series.

\begin{definition}[Laplace transform]
	Let $\theta \in \R$ and suppose that $\psi(\xi)$ is analytic in a half-strip
	\begin{equation}
		\{\xi \in \C \colon \operatorname{dist}(\xi, e^{i\theta} \R_+) < \varepsilon\},
	\end{equation}
	for some $\varepsilon > 0$, and is of exponential type less than $\alpha \in \R$. The \emph{Laplace transform in the direction $\theta$} is defined as
	\begin{equation}
		\mathcal{L}_\theta (\psi) (\tau)
		= \int_{0}^{e^{i\theta} \infty} e^{-\xi/\tau} \psi(\xi)\, d\xi,
	\end{equation}
	whenever the integral converges. In particular, $\mathcal{L}_\theta (\psi) (\tau)$ defines a holomorphic function for
\begin{equation}
	\Re(\tau e^{-i\theta}) > \max(\alpha,0).
\end{equation}
\end{definition}

If the Borel transform of a Gevrey-$1$ series $\psi(\tau)$ and admits analytic continuation with at most exponential growth in the direction $\theta$, then the composition
\begin{equation}
	\S_\theta (\psi) (\tau) \coloneqq (\mathcal{L}_\theta \circ \B) (\psi) (\tau)
\end{equation}
is called the \emph{Borel--Laplace sum} of the series in the direction $\theta$.

The Borel--Laplace sum is compatible with convergent power series. Precisely, under suitable analyticity and growth assumptions, Borel--Laplace summation recovers the original analytic function from its Taylor series, which we prove through the following \namecref{prop:BL-monomial}.

\begin{proposition} \label{prop:BL-monomial}
	The Borel--Laplace sum of the monomials $\tau^{n + 1}$ recovers itself; that is,
	\begin{equation}
		\S_\theta (\tau^{n + 1}) = \tau^{n + 1}.
	\end{equation}
\end{proposition}

\begin{proof}
	By definition of the Borel transform,
	\begin{equation}
		\B(\tau^{n + 1}) = \frac{\xi^n}{n!}.
	\end{equation}
	Taking the Laplace transform then parameterizing the ray $\e^{i\theta} \R_{\ge 0}$ by $\xi = t e^{i\theta}$ for $t \ge 0$ yields
	\begin{equation}
		\mathcal{L}_\theta \left(\frac{\xi^n}{n!}\right) (\tau)
		= \int_{0}^{e^{i\theta} \infty} e^{-\xi/\tau}\, \frac{\xi^n}{n!}\, d\xi
		= \frac{e^{i(n + 1) \theta}}{n!} \int_{0}^{\infty} e^{-t e^{i\theta}/\tau} t^n\, dt,
	\end{equation}
	which, by substituting $u = t e^{i\theta}/\tau$ and $dt = \tau e^{-i\theta}\, du$, becomes
	\begin{equation}
		\mathcal{L}_\theta \left(\frac{\xi^n}{n!}\right) (\tau)
		= \frac{e^{i(n + 1) \theta}}{n!} \cdot \tau^{n + 1} e^{-i(n + 1) \theta} \int_{0}^{\infty} e^{-u} u^n\, du
		= \frac{\tau^{n + 1}}{n!}\, \Gamma(n + 1)
		= \tau^{n + 1}. \tag*{\qedhere}
	\end{equation}
\end{proof}

\begin{corollary} \label{cor:BL-analytic}
	The Borel--Laplace sum of an analytic function recovers the function itself.
\end{corollary}

\begin{proof}
	Follows by applying \zcref{prop:BL-monomial} to the Taylor expansion of the function.
\end{proof}

\paragraph{Median resummation}
A different phenomenon occurs when the Borel transform has singularities along a direction $\theta$. In this case, the Laplace integral along that direction may not be defined, and one instead considers integration along the rays immediately above and below it. Let $\S_{\theta_{+}} (\psi)$ and $\S_{\theta_{-}} (\psi)$ denote the corresponding lateral Borel--Laplace sums, where $\theta_+ = \theta + \varepsilon$ and $\theta_- = \theta - \varepsilon$ for sufficiently small $\varepsilon > 0$. We assume that both lateral sums exist and admit boundary values on a common domain.

\begin{definition}[Stokes discontinuity]
	The \emph{Stokes discontinuity} of $\psi$ across the direction $\theta$ is defined as
	\begin{equation}
		\disc_{\theta} \psi(\tau) = \S_{\theta_{+}} (\psi) (\tau) - \S_{\theta_{-}} (\psi) (\tau),
	\end{equation}
	where the lateral sums are understood as their limiting values as $\varepsilon \to 0^+$.
\end{definition}

The two lateral sums may differ because their integration contours pass on opposite sides of the singularities of the Borel transform. Their difference measures the resulting Stokes phenomenon. When both lateral sums exist on a common domain, their average, known as the \emph{median resummation}, is
\begin{equation} \label{eq:med-resum-disc}
	\Smed_{\theta} (\psi) (\tau)
		\coloneqq \frac{\S_{\theta_+} (\psi) (\tau) + \S_{\theta_-} (\psi) (\tau)}{2}
		=
		\begin{cases}
			\S_{\theta_{+}} (\psi) (\tau) - \frac{1}{2} \disc_{\theta} \psi(\tau), & \Re(e^{i\theta_{+}} \tau) > 0, \\
			\S_{\theta_{-}} (\psi) (\tau) + \frac{1}{2} \disc_{\theta} \psi(\tau), & \Re(e^{i\theta_{-}} \tau) > 0.
		\end{cases}
\end{equation}

\section{Resurgence of \texorpdfstring{$\tf_\r$}{f̃ᵣ}} \label{sec:fr-resurgence}
In this section, we will find the the Borel--Laplace sums of the series $\psi_\r$ and $\psi_{r_2}$, defined in \zcref{lem:fr-asymp}, along $\theta = 0$ and $\theta = \pi$, as well as their median resummation along $\theta = \pi/2$.

\subsection{Borel transforms} \label{sec:borel-transforms}
Our proofs in this section rely on the following formulas.
\begin{itemize}
	\item The Fourier expansion of the Bernoulli polynomials:
	\begin{equation} \label{eq:bernoulli-fourier}
		B_m (x) = -m! \sum_{n \ne 0} \frac{\e(nx)}{(2\pi i n)^m},
		\quad 0 \le x \le 1,\ m \ge 2,
	\end{equation}
	which converges absolutely and uniformly on $[0,1]$.
	
	\item The generating function of the Bernoulli numbers:
	\begin{equation} \label{eq:bernoulli-nums-gen}
		\frac{t}{e^t - 1} = \sum_{k=0}^{\infty} \frac{B_k}{k!}\, t^k, \quad \abs{t} < 2\pi.
	\end{equation}
	
	\item The generating function of $\Li_{-k} (z)$ for $k \ge 0$:
	\begin{equation} \label{eq:polylog-gen}
		-\frac{1}{1 - z^{-1} e^{-t}}
		= \frac{z e^t}{1 - z e^t}
		= \sum_{k=0}^{\infty} \frac{\Li_{-k} (z)}{k!}\, t^k,
		\quad \abs{t} < \min_{n \in \Z} \abs{2\pi i n - \log z}.
	\end{equation}
\end{itemize}
The last identity can be proved by noting that
\begin{equation}
	\frac{z e^t}{1 - z e^t} = \Li_0 (z e^t),
\end{equation}
which, using $w \Li_{-k}' (w) = \Li_{-k - 1} (w)$, implies that
\begin{equation}
	\frac{\pd^k}{\pd t^k} \left(\frac{z e^t}{1 - z e^t}\right)\bigg|_{t = 0}
	= \Li_{-k} (z e^t)\big|_{t=0}
	= \Li_{-k} (z).
\end{equation}
The radius of convergence of the series is the distance to the nearest pole from the origin.

\begin{proposition} \label{prop:borel-transform-psi-r}
	For $\r = \prmatrix*{r_1}{r_2} \in (\R \sm \Z) \times [0,1]$, the series defining $\B(\psi_\r) (\xi)$ converges in a neighborhood of $\xi = 0$ and meromorphically continues to $\C$ via
	\begin{equation}
		\B(\psi_\r) (\xi)
		= \frac{1}{2\pi i} \sum_{n \ne 0} \frac{1}{n^2} \cdot \frac{\e(n r_2)}{1 - \e(r_1) e^{-\xi/n}}.
	\end{equation}
\end{proposition}

\begin{proof}
	By definition, the Borel transform of $\psi_\r$ is
	\begin{equation}
		\B(\psi_\r) (\xi)
		= \sum_{k=0}^{\infty} \frac{(2\pi i)^{k + 1} B_{k + 2} (r_2) \Li_{-k} (\e(-r_1))}{k! (k + 2)!}\, \xi^k.
	\end{equation}
	Using \eqref{eq:bernoulli-fourier}, we can write
	\begin{align}
		\B(\psi_\r) (\xi)
		&= -\sum_{k=0}^{\infty} \frac{(2\pi i)^{k + 1} \Li_{-k} (\e(-r_1))}{k! (k + 2)!}
		\left((k + 2)! \sum_{n \ne 0} \frac{\e(n r_2)}{(2\pi i n)^{k + 2}}\right) \xi^k \\
		&= -\frac{1}{2\pi i} \sum_{k=0}^{\infty} \sum_{n \ne 0} \frac{\e(n r_2) \Li_{-k} (\e(-r_1))}{n^{k + 2} k!}\, \xi^k.
	\end{align}
	Since $\zeta(k + 2) \le \zeta(2)$ for $k \ge 2$, it follows that
	\begin{align}
		\sum_{k=0}^{\infty} \sum_{n \ne 0} \abs*{\frac{\e(n r_2) \Li_{-k} (\e(-r_1))}{n^{k + 2} k!}\, \xi^k}
		&= \sum_{k=0}^{\infty} \frac{2 \zeta(k + 2) \abs{\Li_{-k} (\e(-r_1))}}{k!}\, \abs{\xi}^k \\
		&\le 2\zeta(2) \sum_{k=0}^{\infty} \frac{\abs{\Li_{-k} (\e(-r_1))}}{k!}\, \abs{\xi}^k.
	\end{align}
	The last series has positive radius of convergence by \eqref{eq:polylog-gen}, so $\B(\psi_\r) (\xi)$ converges in a neighborhood of $\xi = 0$. By absolute convergence, we may interchange the sums then use \eqref{eq:polylog-gen} to get
	\begin{equation}
		\B(\psi_\r) (\xi)
		= -\frac{1}{2\pi i} \sum_{n \ne 0} \frac{\e(n r_2)}{n^2} \sum_{k=0}^{\infty} \frac{\Li_{-k} (\e(-r_1))}{k!}\, (\xi/n)^k
		= \frac{1}{2\pi i} \sum_{n \ne 0} \frac{\e(n r_2)}{n^2} \cdot \frac{1}{1 - \e(r_1) e^{-\xi/n}},
	\end{equation}
	which converges absolutely and locally uniformly on $\C \sm \{2\pi i n(k + r_1) \colon n \in \Z_{\ne 0},\ k \in \Z\}$.
\end{proof}

\begin{proposition} \label{prop:borel-transform-psi-r2}
	For $r_2 \in (0,1]$, the series defining $\B(\psi_{r_2}) (\xi)$ converges in a neighborhood of $\xi = 0$ and meromorphically continues to $\C$ via
	\begin{equation}
		\B(\psi_{r_2}) (\xi)
		= \frac{1}{2\pi i} \sum_{n \ne 0} \frac{\e(n r_2)}{n^2} \left(\frac{1}{1 - e^{-\xi/n}} - \frac{n}{\xi}\right).
	\end{equation}
\end{proposition}

\begin{proof}
	The proof of this result is similar to the proof of \zcref{prop:borel-transform-psi-r}. By definition of the Borel transform and \eqref{eq:bernoulli-fourier},
	\begin{equation}
		\B(\psi_{r_2}) (\xi)
		= -\sum_{k=1}^{\infty} \frac{(2\pi i)^k B_k (1) B_{k + 1} (r_2)}{k! (k + 1)!}\, \xi^{k - 1}
		= \frac{1}{2\pi i} \sum_{k=1}^{\infty} \sum_{n \ne 0} \frac{\e(n r_2) B_k (1)}{n^{k + 1} k!}\, \xi^{k - 1}.
	\end{equation}
	This series converges absolutely in a neighborhood of $\xi = 0$ because
	\begin{equation}
		\sum_{k=1}^{\infty} \sum_{n \ne 0} \abs*{\frac{\e(n r_2) B_k (1)}{n^{k + 1} k!}\, \xi^{k - 1}}
		= \sum_{k=1}^{\infty} \frac{2 \zeta(k + 1) \abs{B_k (1)}}{k!}\, \abs{\xi}^{k - 1}
		\le 2 \zeta(2) \sum_{k=1}^{\infty} \frac{\abs{B_k (1)}}{k!}\, \abs{\xi}^{k - 1},
	\end{equation}
	and the latter series has a positive radius of convergence due to \eqref{eq:bernoulli-nums-gen}. Interchanging the sums then using $B_1 (1) = B_1 + 1$ and $B_k (1) = B_k$ for $k \ge 2$ yields
	\begin{equation}
		\B(\psi_{r_2}) (\xi)
		= \frac{1}{2\pi i} \sum_{n \ne 0} \frac{\e(n r_2)}{n^2} \sum_{k=1}^{\infty} \frac{B_k (1)}{k!}\, (\xi/n)^{k - 1}
		= \frac{1}{2\pi i} \sum_{n \ne 0} \frac{\e(n r_2)}{n^2} \left(1 + \sum_{k=1}^{\infty} \frac{B_k}{k!}\, (\xi/n)^{k - 1}\right).
	\end{equation}
	By \eqref{eq:bernoulli-nums-gen}, we have
	\begin{equation}
		\sum_{k=1}^{\infty} \frac{B_k}{k!}\, t^{k - 1}
		= \frac{1}{t} \left(\frac{t}{e^{t} - 1} - B_0\right)
		= \frac{1}{e^t - 1} - \frac{1}{t}
		= \frac{1}{1 - e^{-t}} - \frac{1}{t} - 1.
	\end{equation}
	Hence,
	\begin{equation}
		\B(\psi_{r_2}) (\xi)
		= \frac{1}{2\pi i} \sum_{n \ne 0} \frac{\e(n r_2)}{n^2} \left(\frac{1}{1 - e^{-\xi/n}} - \frac{n}{\xi}\right),
	\end{equation}
	which converges absolutely and locally uniformly on $\C \sm 2\pi i \Z_{\ne 0}$. Note that the singularity at $\xi = 0$ is removable because
	\begin{equation}
		\lim_{\xi \to 0} \left(\frac{1}{1 - e^{-\xi/n}} - \frac{n}{\xi}\right) = 0. \tag*{\qedhere}
	\end{equation}
\end{proof}

\subsection{Resurgent structures} \label{sec:resurge-struct}
\begin{proposition} \label{prop:fr-stokes}
	Let $\r = \prmatrix*{r_1}{r_2} \in (\R \times [0,1]) \sm (\Z \times \{0\})$. The formal power series $\psi_\r$ and $\psi_{r_2}$ in the asymptotic series $\tf_\r$ are Gevrey-$1$ series. Moreover, they are simple resurgent series with simple poles along the imaginary axis at
	\begin{equation}
		\xi_{n,k} = 2\pi i n(k + r_1), \quad n \in \Z_{\ne 0},\ \ k \in \Z, \ \ k + r_1 \ne 0.
	\end{equation}
	If $y = n(k + r_1) = \xi_{n,k} / (2\pi i)$, then the corresponding Stokes constant is
	\begin{equation}
		S_y = -\sum_{\scalewidth{\substack{n \in \Z_{\ne 0} \\ y / n - r_1 \in \Z}}} \frac{\e(n r_2)}{n}.
	\end{equation}
\end{proposition}

\begin{proof}
	We have shown in the proofs of \zcref{prop:borel-transform-psi-r,prop:borel-transform-psi-r2} that series defining $\B(\psi_\r) (\xi)$ and $\B(\psi_{r_2}) (\xi)$ have positive radii of convergence. This implies (e.g.~by the divergence test for infinite series) that $\psi_\r$ and $\psi_{r_2}$ are Gevrey-$1$ series.
	
	Suppose now that $r_1 \notin \Z$. Then the series appearing in $\tf_\r$ is $\psi_\r$. By \zcref{prop:borel-transform-psi-r},
	\begin{equation}
		\B(\psi_\r) (\xi)
		= \frac{1}{2\pi i} \sum_{n \ne 0} \frac{1}{n^2} \cdot \frac{\e(n r_2)}{1 - \e(r_1) e^{-\xi/n}}.
	\end{equation}
	For fixed $n \in \Z_{\ne 0}$, the summand has simple poles at
	\begin{equation}
		\xi_{n,k} = 2\pi i n (k + r_1), \quad k \in \Z,
	\end{equation}
	with residue $\e(n r_2) / n$. But it is possible (when $r_1 \in \Q$) for different terms of the sum have a pole at the same point, so we must collect the residues of all such terms to obtain the residue of $\B(\psi_\r) (\xi)$ at $\xi = \xi_{n,k}$. Write $\xi_{n,k} = 2\pi i y$ with $y = n(k + r_1)$. Then the corresponding Stokes constant is
	\begin{equation}
		S_y
		= -2\pi i \Res_{\xi = \xi_{n,k}} \B(\psi_\r) (\xi)
		= -2\pi i \cdot \frac{1}{2\pi i} \sum_{y = n(k + r_1)} \frac{\e(n r_2)}{n}
		= -\sum_{\scalewidth{y = n(k + r_1)}} \frac{\e(n r_2)}{n},
	\end{equation}
	where the sum runs over all pairs $(n,k)$ with $n \in \Z_{\ne 0}$ and $k \in \Z$ such that $y = n(k + r_1)$. By writing $k = y/n - r_1$, we rewrite the sum as
	\begin{equation}
		S_y = -\sum_{\scalewidth{\substack{n \in \Z_{\ne 0} \\ y / n - r_1 \in \Z}}} \frac{\e(n r_2)}{n}.
	\end{equation}
	
	Now let $r_1 \in \Z$. Then the series appearing in $\tf_\r$ is $\psi_{\r_2}$. By \zcref{prop:borel-transform-psi-r2},
	\begin{equation}
		\B(\psi_{r_2}) (\xi)
		= \frac{1}{2\pi i} \sum_{n \ne 0} \frac{\e(n r_2)}{n^2} \left(\frac{1}{1 - e^{-\xi/n}} - \frac{n}{\xi}\right).
	\end{equation}
	Here, each summand has simple poles at
	\begin{equation}
		\xi_{n,k} = 2\pi i nk, \quad k \in \Z_{\ne 0}
	\end{equation}
	and a removable singularity at $\xi = 0$. Note that since $r_1 \in \Z$, this is equivalent to
	\begin{equation}
		\xi_{n,k} = 2\pi i n (k + r_1), \quad k \in \Z, \ \ k + r_1 \ne 0.
	\end{equation}
	The residue of $n$-th summand at $\xi = \xi_{n,k}$ is again $\e(n r_2) / n$, so the Stokes constant at $\xi = 2\pi i y$ for $y \ne 0$ is again equal to
	\begin{equation}
		S_y
		= -\sum_{\scalewidth{\substack{n \in \Z_{\ne 0} \\ y / n \in \Z}}} \frac{\e(n r_2)}{n}
		= -\sum_{\scalewidth{\substack{n \in \Z_{\ne 0} \\ y / n - r_1 \in \Z}}} \frac{\e(n r_2)}{n}. \tag*{\qedhere}
	\end{equation}
\end{proof}

Following \eqref{eq:L-pm-general}, we consider the pair of general Dirichlet series indexed by the set of distinct $y > 0$ such that
\begin{equation}
	y = n (k + r_1), \quad n \in \Z_{\ne 0}, \ \ k \in \Z, \ \ k + r_1 \ne 0,
\end{equation}
with coefficients the Stokes constants $S_{\pm y}$ from \zcref{prop:fr-stokes}:
\begin{equation} \label{eq:fr-L-series}
	L_+ (s) \coloneqq \sum_{y > 0} \frac{S_y}{y^s}
	\qqand
	L_- (s) \coloneqq -\sum_{y > 0} \frac{S_{-y}}{y^s}.
\end{equation}

\begin{proposition}
	The Dirichlet series $L_\pm (s)$ in \eqref{eq:fr-L-series} converge absolutely for $\Re(s) > 1$. Moreover, both series have a meromorphic continuation to the complex plane via
	\begin{align}
		L_+ (s) &= \zeta(s, 1 - \{r_1\}) \Li_{s + 1} (\e(-r_2)) - \zeta(s, 1 - \{-r_1\}) \Li_{s + 1} (\e(r_2)), \\
		L_- (s) &= \zeta(s, 1 - \{r_1\}) \Li_{s + 1} (\e(r_2)) - \zeta(s, 1 - \{-r_1\}) \Li_{s + 1} (\e(-r_2)),
	\end{align}
	where $\zeta(s,\alpha)$ is the Hurwitz zeta function and $\{x\} \coloneqq x - \floor{x}$ denotes the fractional part of $x$.
\end{proposition}

\begin{proof}
	Let $\rho = \Re(s) > 1$. Observe that
	\begin{equation}
		\abs{S_y}
		\le \sum_{\scalewidth[0.7]{\substack{n \in \Z_{\ne 0} \\ y / n - r_1 \in \Z}}} \abs*{\frac{\e(n r_2)}{n}}
		= \sum_{\scalewidth[0.7]{\substack{n \in \Z_{\ne 0} \\ y / n - \{r_1\} \in \Z}}} \frac{1}{\abs{n}}.
	\end{equation}
	If $r_1 \in \R \sm \Q$, then there is at most one $n \in \Z_{\ne 0}$ such that $y/n - \{r_1\} \in \Z$ (indeed, if $y/n_j - \{r_1\} \in \Z$ for $j = 1,2$, then $(n_1 - n_2) \{r_1\} \in \Z$ and $n_1 = n_2$), so $\abs{S_y} \le 1$. Therefore,
	\begin{align}
		\abs{L_\pm (s)}
		\le \sum_{y > 0} \abs*{\frac{S_{\pm y}}{y^s}}
		= \sum_{y > 0} \frac{1}{\abs{y}^\rho}
		&\le \sum_{n=1}^{\infty} \sum_{k=0}^{\infty} \frac{1}{(n (k + \{r_1\}))^\rho}
		+ \sum_{n=1}^{\infty} \sum_{k=0}^{\infty} \frac{1}{(n (k + 1 - \{r_1\}))^\rho} \\
		&= \zeta(\rho) (\zeta(\rho, \{r_1\}) + \zeta(\rho, 1 - \{r_1\}))
		< \infty.
	\end{align}
	
	Now let $r_1 \in \Q$ and write $r_1 = p/q$ in lowest terms. Then the condition on the sum becomes
	\begin{equation}
		\frac{y}{n} - \frac{p}{q} = k
	\end{equation}
	for some $k \in \Z$. Rearranging yields $qy = n (qk + p)$, so the sum runs over a subset of all $n \in \Z_{\ne 0}$ such that $n \mid qy$. Therefore, with $m = \abs{qy}$,
	\begin{equation}
		\abs{S_y}
		\le \sum_{\substack{n \in \Z_{\ne 0} \\ n \mid qy}} \frac{1}{\abs{n}}
		= \frac{1}{\abs{qy}} \sum_{\substack{d \in \Z_{\ne 0} \\ d \mid qy}} \abs{d}
		= \frac{2 \sigma_1 (m)}{m},
	\end{equation}
	where $\sigma_1 (m) \coloneqq \sum_{d \mid m} d$ is the sum of divisors function. Consequently,
	\begin{equation}
		\abs{L_\pm (s)}
		\le \sum_{y > 0} \abs*{\frac{S_{\pm y}}{y^s}}
		= 2\sum_{m = 1} \frac{\sigma_1 (m) / m}{(m / q)^\rho}
		= 2q^\rho \sum_{m=1}^{\infty} \sum_{d \mid m} \frac{d}{m^{\rho + 1}}.
	\end{equation}
	Substituting $m = jd$ yields
	\begin{equation}
		\abs{L_\pm (s)}
		\le 2q^\rho \sum_{j=1}^{\infty} \sum_{d=1}^{\infty} \frac{d}{(jd)^{\rho + 1}}
		= 2q^\rho \sum_{j=1}^{\infty} \frac{1}{j^{\rho + 1}} \sum_{d=1}^{\infty} \frac{1}{d^\rho}
		= 2q^\rho \zeta(\rho) \zeta(\rho + 1) < \infty.
	\end{equation}
	This proves that both series converge absolutely for $\Re(s) > 1$.
	
	Recall that the sums in \eqref{eq:fr-L-series} run over distinct $y > 0$ such that
	\begin{equation}
		y = n (k + r_1), \quad n \in \Z_{\ne 0}, \ \ k \in \Z, \ \ k + r_1 \ne 0.
	\end{equation}
	Therefore,
	\begin{equation}
		L_+ (s)
		= \sum_{y > 0} \frac{S_y}{y^s}
		= -\sum_{n \in \Z_{\ne 0}} \frac{\e(n r_2)}{n} \sum_{\scalewidth{\substack{k \in \Z \\ n(k + r_1) > 0}}} \frac{1}{(n (k + r_1))^s}.
	\end{equation}
	By splitting the sum according to the sign of $n$ then replacing $n$ by $-n$ in the negative half, we obtain
	\begin{align}
		L_+ (s)
		&= \sum_{n=1}^{\infty} \frac{\e(-n r_2)}{n} \sum_{\scalewidth{\substack{k \in \Z \\ k + r_1 < 0}}} \frac{1}{(-n (k + r_1))^s}
		- \sum_{n=1}^{\infty} \frac{\e(n r_2)}{n} \sum_{\scalewidth{\substack{k \in \Z \\ k + r_1 > 0}}} \frac{1}{(n (k + r_1))^s} \\
		&= \sum_{n=1}^{\infty} \frac{\e(-n r_2)}{n^{s + 1}} \sum_{\scalewidth{\substack{k \in \Z \\ k + r_1 < 0}}} \frac{1}{(-k - r_1)^s}
		- \sum_{n=1}^{\infty} \frac{\e(n r_2)}{n^{s + 1}} \sum_{\scalewidth{\substack{k \in \Z \\ k + r_1 > 0}}} \frac{1}{(k + r_1)^s} \\
		&= \Li_{s + 1} (\e(-r_2)) \sum_{\scalewidth{\substack{k \in \Z \\ k - r_1 > 0}}} \frac{1}{(k - r_1)^s}
		- \Li_{s + 1} (\e(r_2)) \sum_{\scalewidth{\substack{k \in \Z \\ k + r_1 > 0}}} \frac{1}{(k + r_1)^s},
	\end{align}
	where in the last step we replaced $k$ by $-k$ in the first sum and also used
	\begin{equation}
		\Li_{s} (z) = \sum_{n=1}^{\infty} \frac{z^n}{n^s}.
	\end{equation}
	The conditions $k - r_1 > 0$ and $k + r_1 > 0$ are equivalent to $k \ge \floor{r_1} + 1$ and $k > 1 - \ceil{r_1}$ respectively. Thus, using $x - \floor{x} = \{x\}$ and $x - \ceil{x} + 1 = 1 - \{-x\}$, we find
	\begin{align}
		L_+ (s)
		&= \Li_{s + 1} (\e(-r_2)) \sum_{k = \floor{r_1} + 1}^{\infty} \frac{1}{(k - r_1)^s}
		- \Li_{s + 1} (\e(r_2)) \sum_{k = 1 - \ceil{r_1}}^{\infty} \frac{1}{(k + r_1)^s} \\
		&= \Li_{s + 1} (\e(-r_2)) \sum_{k=0}^{\infty} \frac{1}{(k + 1 - \{r_1\})^s}
		- \Li_{s + 1} (\e(r_2)) \sum_{k=0}^{\infty} \frac{1}{(k + 1 - \{-r_1\})^s} \\
		&= \zeta(s, 1 - \{r_1\}) \Li_{s + 1} (\e(-r_2)) - \zeta(s, 1 - \{-r_1\}) \Li_{s + 1} (\e(r_2)).
	\end{align}
	An analogous argument proves the corresponding identity for $L_- (s)$.
\end{proof}

\subsection{Borel--Laplace sums} \label{sec:BL-sums}
Before we calculate the Borel--Laplace sums of $\tf_\r$, we prove an integral identity that will be used repeatedly in the computations. Its proof relies on evaluating a suitable contour integral.
\begin{lemma} \label{lem:int-cont-shift}
	Let $\r = \prmatrix*{r_1}{r_2} \in (\R \sm \Z) \times [0,1)$ and $\tau \in \H$ such that $0 < \Re(\tau) < \{r_1\} / r_2$ \textup{(}with $\{r_1\} / r_2 = \infty$ if $r_2 = 0$\textup{)}. Then
	\begin{equation}
		\int_{-\infty}^{\infty} \frac{\ell(t/\tau + 2\pi i r_2)}{1 - \e(r_1) e^t}\, dt
		= \int_{-\infty}^{\infty} \frac{\ell(t/\tau)}{1 - \e(r_1 - r_2 \tau) e^t}\, dt,
	\end{equation}
	where $\ell(z)$ is the logarithm defined in \zcref{prop:continuous-L}.
\end{lemma}

\begin{proof}
	If $r_2 = 0$, then there is nothing to prove, so suppose that $r_2 \ne 0$. First, we replace $t$ by $-t + 2\pi i r_2 \tau$ in the integral on the left-hand side to get
	\begin{equation}
		\int_{-\infty}^{\infty} \frac{\ell(t/\tau + 2\pi i r_2)}{1 - \e(r_1) e^t}\, dt
		= \int_{-\infty + ih}^{\infty + ih} \frac{\ell(t/\tau)}{1 - \e(r_1 - r_2 \tau) e^t}\, dt,
	\end{equation}
	where $h = 2\pi r_2 \Re(\tau)$. Now we need to shift the contour of integration to the real line.
	
	Consider the integral
	\begin{equation}
		\int_{C_{T, \varepsilon}} \frac{\ell(t/\tau)}{1 - \e(r_1 - r_2 \tau) e^t}\, dt,
	\end{equation}
	where $T > 0$, $\varepsilon > 0$, and $C_{T, \varepsilon}$ is the rectangular contour with vertices $\pm T$ and $\pm T + ih$ and a semicircular perturbation of radius $\varepsilon$ above the origin, taken to avoid the branch point at the origin (see \zcref{fig:int-contour}). To prove the result, we will use the residue theorem to show that the integral vanishes for any $T > 0$ and $\varepsilon > 0$, then we will show that the integrals on $C_{T}^{(2)}$, $C_{T}^{(4)}$, and $C_\varepsilon$ vanish as $T \to \infty$ and $\varepsilon \to 0^+$.
	
	To apply the residue theorem, we need to show that $\ell(t/\tau)$ is meromorphic inside the contour $C_{T,\varepsilon}$ and analytic on the contour. For meromorphy, it suffices to show that if $t$ is inside the contour and $t/\tau = iy$ for some $y \in \R$, then $y \in (0,2\pi)$. Indeed, if that is the case, then $t = iy \tau$ and $\Im(t) = y \Re(\tau)$. But we have that $0 \le \Im(t) \le h = 2\pi r_2 \Re(\tau)$, so
	\begin{equation}
		0 \le y \le 2\pi r_2 < 2\pi.
	\end{equation}
	If $y = 0$, then $t = 0$, which is not possible due to the perturbation near zero on the contour. Hence $y \in (0,2\pi)$.
	
	Now we are ready to apply the residue theorem. The integrand has poles whenever $e^t = \e(r_2 \tau - r_1)$, or equivalently when $t = t_n = 2\pi i (r_2 \tau - r_1 + n)$ for $n \in \Z$. A pole at $t = t_n$ occurs in or on $C_{T, \varepsilon}$ if and only if
	\begin{equation}
		r_1 - r_2 \Re(\tau) \le n \le r_1.
	\end{equation}
	Since $r_1 \notin \Z$ and $0 < \Re(\tau) < \{r_1\} / r_2$, no such integer $n$ exists. Thus, by the residue theorem
	\begin{equation} \label{eq:F-r-tau-contour-int}
		\int_{C_{T, \varepsilon}} F(t)\, dt = 0,
		\quad \text{where} \quad
		F(t) = \frac{\ell(t/\tau)}{1 - \e(r_1 - r_2 \tau) e^t}.
	\end{equation}
	
	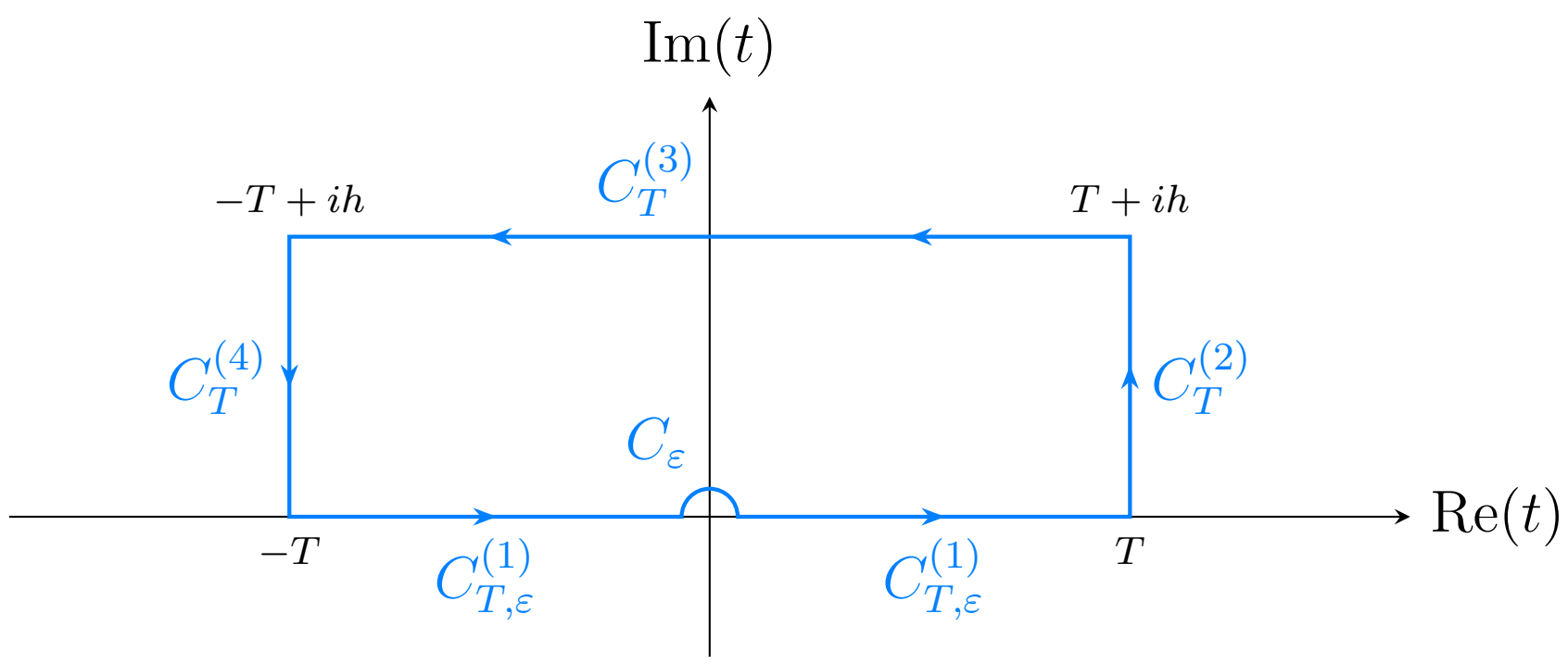
\begin{figure}
		\centering
		\begin{tikzpicture}
			\tikzmath{
				\x1 = -5;
				\x2 = 5;
				\y1 = -1;
				\y2 = 3;
				\s = 0;
			}
			
			\tikzset{
				mid arrow/.style={
					sloped, inner sep=0pt, pos=#1,
					node contents={\tikz\draw[-{Stealth}] (0,0) -- (0.01,0);}
				},
				mid arrow/.default=0.5,
				mid arrow rev/.style={
					sloped, inner sep=0pt, pos=#1,
					node contents={\tikz\draw[{Stealth}-] (0,0) -- (0.01,0);}
				},
				mid arrow rev/.default=0.5
			}
			
			\draw[-stealth] (\x1 - \s, 0) -- (\x2 + \s, 0) node[right] {$\Re(t)$};
			\draw[-stealth] (0, \y1 - \s) -- (0, \y2 + \s) node[above] {$\Im(t)$};
			
			\node[below] at (-3,0) {\scriptsize $-T$};
			\node[below] at (3,0) {\scriptsize $T$};
			\node[above] at (3,2) {\scriptsize $T + ih$};
			\node[above] at (-3,2) {\scriptsize $-T + ih$};
			
			\begin{scope}
				\clip (\x1, \y1) rectangle (\x2, \y2);
				
				\draw[thick, blue!50!cyan] (-3,0) -- (-0.2,0) node[mid arrow] node[midway, below] {$C_{T, \varepsilon}^{(1)}$}
				arc[start angle=180, end angle=0, radius=0.2] node[midway, above left] {$C_\varepsilon$}
				-- (3,0) node[mid arrow] node[midway, below] {$C_{T, \varepsilon}^{(1)}$}
				-- (3,2) node[mid arrow] node[midway, right] {$C_{T}^{(2)}$}
				-- (-3,2) node[mid arrow rev = 0.25] node[mid arrow rev = 0.75] node[pos=0.575, above] {$C_{T}^{(3)}$}
				-- cycle node[mid arrow] node[midway, left] {$C_{T}^{(4)}$};
			\end{scope}
		\end{tikzpicture}
		\caption{The contour $C_{T,\varepsilon}$.}
		\label{fig:int-contour}
	\end{figure}
	
	Now we consider the integrals on $C_\varepsilon$, $C_T^{(2)}$, and $C_T^{(4)}$ as $T \to \infty$ and $\varepsilon \to 0^+$. If $t$ is on $C_T^{(2)}$, then, as $T \to \infty$,
	\begin{equation}
		\ell(t/\tau)
		= \frac{t}{\tau} - \pi i + \log(1 - e^{-t/\tau})
		= \frac{t}{\tau} + O(1)
		= O(T)
	\end{equation}
	and
	\begin{equation}
		\abs{1 - \e(r_1 - r_2 \tau) e^t} \asymp e^T,
	\end{equation}
	so $F(t) = O(T e^{-T})$. It follows that
	\begin{equation} \label{eq:int-CT2}
		\lim_{T \to \infty} \int_{C_T^{(2)}} F(t)\, dt = 0.
	\end{equation}
	For $t$ on $C_T^{(4)}$, there exists $c > 0$ such that $e^{t/\tau} = O(e^{-cT})$ as $T \to \infty$. Thus, for $T$ large enough, we have
	\begin{equation}
		\ell(t/\tau)
		= \log(1 - e^{t/\tau})
		= -e^{t/\tau} + O(e^{2t/\tau})
		= O(e^{-cT}).
	\end{equation}
	Since the denominator of $F(t)$ converges to $1$, it follows that $F(t) = O(e^{-cT})$, which implies that
	\begin{equation} \label{eq:int-CT4}
		\lim_{T \to \infty} \int_{C_T^{(4)}} F(t)\, dt = 0.
	\end{equation}
	
	Now consider $t$ on $C_\varepsilon$. Since
	\begin{equation}
		1 - e^{\pm t/\tau}
		= \mp\frac{t}{\tau} + O(t^2)
		= O(\varepsilon)
		\quad \text{as } \varepsilon \to 0^+,
	\end{equation}
	it follows that
	\begin{equation}
		\ell(t/\tau)
		= O(\log \varepsilon) + O(1)
		= O(\log \varepsilon)
		\quad \text{as } \varepsilon \to 0^+.
	\end{equation}
	Furthermore, since $\Re(\tau) < \{r_1\} / r_2$, there exists $\delta > 0$, not depending on $t$, such that
	\begin{equation}
		\abs{1 - \e(r_1 - r_2 \tau) e^t} > \delta
	\end{equation}
	Hence $F(t) = O(\log \varepsilon)$ and
	\begin{equation} \label{eq:int-C-epsilon}
		\lim_{\varepsilon \to 0^+} \abs*{\int_{C_\varepsilon} F(t)\, dt}
		\le \lim_{\varepsilon \to 0^+} \operatorname{length}(C_\varepsilon) \cdot O(\log \varepsilon)
		= \lim_{\varepsilon \to 0^+} \pi \varepsilon \cdot O(\log \varepsilon)
		= 0.
	\end{equation}
	
	By sending $T \to \infty$ and $\varepsilon \to 0^+$ in \eqref{eq:F-r-tau-contour-int} and using \eqref{eq:int-CT2}, \eqref{eq:int-CT4}, and \eqref{eq:int-C-epsilon}, we obtain
	\begin{equation}
		\int_{-\infty}^{\infty} F(t)\, dt
		= \lim_{T \to \infty} \lim_{\varepsilon \to 0^+} \int_{C_{T, \varepsilon}^{(1)}} F(t)\, dt
		= \lim_{T \to \infty} \int_{-C_T^{(3)}} F(t)\, dt
		= \int_{-\infty + ih}^{\infty + ih} F(t)\, dt.
	\end{equation}
	Recall that
	\begin{equation}
		F(t)
		= \frac{\ell(t/\tau)}{1 - \e(r_1 - r_2 \tau) e^t},
	\end{equation}
	so
	\begin{equation}
		\int_{-\infty + ih}^{\infty + ih} \frac{\ell(t/\tau)}{1 - \e(r_1 - r_2 \tau) e^t}\, dt
		= \int_{-\infty}^{\infty} \frac{\ell(t/\tau)}{1 - \e(r_1 - r_2 \tau) e^t}\, dt. \tag*{\qedhere}
	\end{equation}
\end{proof}

\begin{remark} \label{rmk:int-cont-shift}
	If $\r \in \Z \times [0,1)$ in \zcref{lem:int-cont-shift}, then the integral
	\begin{equation}
		\int_{-\infty}^{\infty} \frac{\ell(t/\tau + 2\pi i r_2)}{1 - \e(r_1) e^t}\, dt
		= \int_{-\infty}^{\infty} \frac{\ell(t/\tau + 2\pi i r_2)}{1 - e^t}\, dt
	\end{equation}
	diverges because the integrand has a pole at $t = 0$. If we instead replace this integral with its principal value, then we get the relation
	\begin{equation}
		\PV \int_{-\infty}^{\infty} \frac{\ell(t/\tau + 2\pi i r_2)}{1 - e^t}\, dt
		= \int_{-\infty}^{\infty} \frac{\ell(t/\tau)}{1 - \e(-r_2 \tau) e^t}\, dt + \pi i \log(1 - \e(r_2)).
	\end{equation}
\end{remark}

\begin{remark}
	The upper bound on $\Re(\tau)$ in \zcref{lem:int-cont-shift} is for mere convenience; its purpose is to ensure that the contour of integration contains no poles of the integrand. If we allow $\Re(\tau) > 0$ and arbitrary $r_1 \in \R$, we obtain the formula
	\begin{equation}
		\PV \int_{-\infty}^{\infty} \frac{\ell(t/\tau + 2\pi i r_2)}{1 - \e(r_1) e^t}\, dt
		= \PV \int_{-\infty}^{\infty} \frac{\ell(t/\tau)}{1 - \e(r_1 - r_2 \tau) e^t}\, dt
		+ 2\pi i R_\r (\tau),
	\end{equation}
	where $R_\r (\tau)$ is the sum of the residues and is equal to
	\begin{align}
		R_\r (\tau)
		= \psum_{\scalewidth[0.7]{-r_1 \le n \le r_2 \Re(\tau) - r_1}} \log(1 - \e(r_2 - (r_1 + n) / \tau))
	\end{align}
	with the prime on the sum indicating that if $r_1 \in \Z$ or $r_2 \Re(\tau) - r_1 \in \Z$, then the corresponding boundary term is to be weighed by $1/2$. Here, the principal values on the integrals are needed only if $r_1 \in \Z$ for the left side and if $r_2 \Re(\tau) - r_1 \in \Z$ for the right side, in which case the respective integrand has a pole on the real line.
\end{remark}

\subsubsection{Borel--Laplace sums for \texorpdfstring{$r_1 \notin \Z$}{r₁ ∉ ℤ}}
\begin{theorem} \label{thm:bl-fr-r1-nonint}
	Let $\r = \prmatrix*{r_1}{r_2} \in (\R \sm \Z) \times \R_{\ge 0}$ and $\tau \in \H$. The Borel--Laplace sums of the asymptotic series of $f_\r (\tau)$ as $\tau \to 0$ along the rays $\theta = 0$ and $\theta = \pi$ satisfy
	\begin{alignat}{2}
		\S_0   (\tf_\r) (\tau) &= f_\r (\tau) - f_{\{ S\r\}} (S \cdot \tau), & \qquad &\Re(\tau) > 0, \label{eq:s0-fr-r1-nonint} \\
		\S_\pi (\tf_\r) (\tau) &= f_\r (\tau) - f_{\{-S\r\}} (S \cdot \tau), & \qquad &\Re(\tau) < 0. \label{eq:s-pi-fr-r1-nonint}
	\end{alignat}
\end{theorem}

\begin{proof}
	We begin by proving \eqref{eq:s0-fr-r1-nonint} then use it to deduce \eqref{eq:s-pi-fr-r1-nonint}. We proceed by induction on $\floor{r_2}$. Let $r_2 \in [0,1)$ (i.e.~$\floor{r_2} = 0$) and suppose that $0 < \Re(\tau) < \{r_1\} / r_2$ (with $\{r_1\} / r_2 = \infty$ if $r_2 = 0$). By \zcref{prop:borel-transform-psi-r}, we have
	\begin{equation}
		\S_0 (\psi_\r) (\tau)
		= \int_{0}^{\infty} e^{-\xi/\tau} \B(\psi_\r) (\xi)\, d\xi
		= \frac{1}{2\pi i} \int_{0}^{\infty} \sum_{n \ne 0} \frac{e^{-\xi/\tau}}{n^2} \cdot \frac{\e(n r_2)}{1 - \e(r_1) e^{-\xi/n}}\, d\xi
	\end{equation}
	For $n \notin \Z_{\ne 0}$ and $\xi \ge 0$, there exists $K > 0$ such that $\abs{1 - \e(r_1) e^{-\xi/n}} \ge K$. Therefore,
	\begin{equation}
		\int_{0}^{\infty} \abs[\bigg]{\sum_{n \ne 0} \frac{e^{-\xi/\tau}}{n^2} \cdot \frac{\e(n r_2)}{1 - \e(r_1) e^{-\xi/n}}}\, d\xi
		\le \int_{0}^{\infty} \sum_{n \ne 0} \frac{e^{-\xi \Re(1/\tau)}}{K n^2}\, d\xi
		= \frac{2\zeta(2)}{K} \int_{0}^{\infty} e^{-\xi \Re(1/\tau)}\, d\xi
		< \infty.
	\end{equation}
	So we may interchange the sum and the integral by Fubini's Theorem to obtain
	\begin{equation}
		\S_0 (\psi_\r) (\tau)
		= \frac{1}{2\pi i} \sum_{n=1}^{\infty} \frac{1}{n^2} \left(
		\int_{0}^{\infty}\frac{e^{-\xi/\tau} \e(n r_2)}{1 - \e(r_1) e^{-\xi/n}}\, d\xi
		+ \int_{0}^{\infty} \frac{e^{-\xi/\tau} \e(-n r_2)}{1 - \e(r_1) e^{\xi/n}}\, d\xi\right).
	\end{equation}
	By substituting $\xi = -nt$ in the first integral and $\xi = nt$ in the second then using the Taylor expansion of $\log(1 - z)$, the above becomes
	\begin{align}
		\S_0 (\psi_\r) (\tau)
		&= \frac{1}{2\pi i} \sum_{n=1}^{\infty} \frac{1}{n}
		\left(\int_{-\infty}^{0} \frac{e^{nt/\tau} \e(n r_2)}{1 - \e(r_1) e^t}\, dt
		+ \int_{0}^{\infty} \frac{e^{-nt/\tau} \e(-n r_2)}{1 - \e(r_1) e^t}\, dt\right) \\
		&= \frac{1}{2\pi i}
		\left(\int_{-\infty}^{0} \frac{1}{1 - \e(r_1) e^t} \sum_{n=1}^{\infty} \frac{e^{nt/\tau} \e(n r_2)}{n}\, dt
		+ \int_{0}^{\infty} \frac{1}{1 - \e(r_1) e^t} \sum_{n=1}^{\infty} \frac{e^{-nt/\tau} \e(-n r_2)}{n}\, dt\right) \\[-0.5\baselineskip]
		&= -\frac{1}{2\pi i} \int_{-\infty}^{0} \frac{\log(1 - e^{t/\tau} \e(r_2))}{1 - \e(r_1) e^t}\, dt
		- \frac{1}{2\pi i} \int_{0}^{\infty} \frac{\log(1 - e^{-t/\tau} \e(-r_2))}{1 - \e(r_1) e^t}\, dt, \label{eq:s0-I1}
	\end{align}
	Since
	\begin{equation}
		\Re\left(\frac{t}{\tau} + 2\pi i r_2\right)
		= \frac{t}{\abs{\tau}^2} \Re(\tau)
	\end{equation}
	and $\Re(\tau) > 0$, it follows that
	\begin{equation} \label{eq:l-reflection}
		\ell\left(\frac{t}{\tau} + 2\pi i r_2\right)
		=
		\begin{cases}
			\log(1 - e^{t/\tau} \e(r_2)),                               & t \le 0, \\
			t/\tau + 2\pi i B_1 (r_2) + \log(1 - e^{-t/\tau} \e(-r_2)), & t > 0,
		\end{cases}
	\end{equation}
	where $B_1 (x) = x - 1/2$. Therefore, we may rewrite \eqref{eq:s0-I1} as
	\begin{align}
		\hspace{-25pt}
		\S_0 (\psi_\r) (\tau)
		&= -\frac{1}{2\pi i} \int_{-\infty}^{0} \frac{\ell(t/\tau + 2\pi i r_2)}{1 - \e(r_1) e^t}\, dt
		- \frac{1}{2\pi i} \int_{0}^{\infty} \frac{\ell(t/\tau + 2\pi i r_2) - t/\tau - 2\pi i B_1 (r_2)}{1 - \e(r_1) e^t}\, dt \\
		&\hspace{-10pt}= -\frac{1}{2\pi i} \int_{-\infty}^{\infty} \frac{\ell(t/\tau + 2\pi i r_2)}{1 - \e(r_1) e^t}\, dt
		+ \frac{1}{2\pi i \tau} \int_{0}^{\infty} \frac{t}{1 - \e(r_1) e^t}\, dt
		+ \int_{0}^{\infty} \frac{B_1 (r_2)}{1 - \e(r_1) e^t}\, dt. \label{eq:s0-I2}
	\end{align}
	Using the formula $\int_{0}^{\infty} \frac{\xi d\xi}{1 - z e^\xi} = -\Li_2 (z^{-1})$, the second integral evaluates to
	\begin{equation}
		\frac{1}{2\pi i \tau} \int_{0}^{\infty} \frac{t}{1 - \e(r_1) e^t}\, dt
		= -\frac{1}{2\pi i \tau} \Li_2 (\e(-r_1)),
	\end{equation}
	and the last integral is
	\begin{equation}
		\int_{0}^{\infty} \frac{B_1 (r_2)}{1 - \e(r_1) e^t}\, dt
		= -B_1 (r_2) \int_{0}^{\infty} \frac{\e(-r_1) e^{-t}}{1 - \e(-r_1) e^{-t}}\, dt
		= B_1 (r_2) \log(1 - \e(-r_1)).
	\end{equation}
	By substituting the above into \eqref{eq:s0-I2} then using \zcref{lem:int-cont-shift}, we get
	\begin{equation} \label{eq:s0-psi}
		\S_0 (\psi_\r) (\tau)
		= -\frac{1}{2\pi i} \int_{-\infty}^{\infty} \frac{\ell(t/\tau)}{1 - \e(r_1 - r_2 \tau) e^t}\, dt
		- \frac{1}{2\pi i \tau} \Li_2 (\e(-r_1))
		+ B_1 (r_2) \log(1 - \e(-r_1)).
	\end{equation}
	Now \zcref{lem:fr-asymp,cor:BL-analytic} yield
	\begin{align}
		\S_0 (\tf_\r) (\tau)
		&= \frac{1}{2\pi i \tau} \Li_2 (\e(-r_1)) - B_1 (r_2) \log(1 - \e(-r_1)) + \S_0 (\psi_\r) (\tau) \\
		&= -\frac{1}{2\pi i} \int_{-\infty}^{\infty} \frac{\ell(t/\tau)}{1 - \e(r_1 - r_2 \tau) e^t}\, dt.
	\end{align}
	Applying \zcref{cor:int-to-qpoch} with $z = r_2 \tau - \{r_1\}$ and noting that $\e(\{r_1\}) = \e(r_1)$ yields
	\begin{equation}
		\S_0 (\tf_\r) (\tau)
		= \slog\frac{\qpoch{\e(r_2 \tau - \{r_1\})}{\e(\tau)}}{\qpoch{\e(r_2 - \{r_1\} / \tau)}{\e(-1/\tau)}}
		= f_\r (\tau) - f_{\{S\r\}} (S \cdot \tau).
	\end{equation}
	Since both sides of the above equation are analytic for $\Re(\tau) > 0$, \eqref{eq:s0-fr-r1-nonint} follows for $r_2 \in [0,1)$.
	
	Now suppose that the result holds for some $\r = \prmatrix*{r_1}{r_2}$ with $r_2 \ge 0$. We will show that it also holds for $\tr \coloneqq \prmatrix*{r_1}{r_2 + 1}$. By definition of $f_\r$, we have
	\begin{equation}
		f_{\tr} (\tau)
		= f_\r (\tau) - \log(1 - \e(r_2 \tau - r_1)).
	\end{equation}
	Therefore, by hypothesis and \zcref{cor:BL-analytic}, taking the Borel--Laplace sum of both sides yields
	\begin{equation}
		\S_0 (\tf_{\tr}) (\tau)
		= \S_0 (\tf_\r) (\tau)
		= f_\r (\tau) - f_{\{S\r\}} (S \cdot \tau) + \log(1 - \e(r_2 \tau - r_1))
		= f_{\tr} (\tau) - f_{\{S\tr\}} (S \cdot \tau)
	\end{equation}
	where we used $\{S\tr\} = \{S\r\}$. This completes the proof of \eqref{eq:s0-fr-r1-nonint}.
	
	Now let $\Re(\tau) < 0$. To prove \eqref{eq:s-pi-fr-r1-nonint}, we use \eqref{eq:s0-fr-r1-nonint} with $\tau \mapsto -\bar{\tau}$ and $\r \mapsto \r' \coloneqq \prmatrix*{-r_1}{r_2}$ to get
	\begin{equation} \label{eq:s0-tau-bar}
		\S_0 (\tf_{\r'}) (-\bar{\tau}) = f_{\r'} (-\bar{\tau}) - f_{\{S\r'\}} (S \cdot (-\bar{\tau})).
	\end{equation}
	To get the result, we take the complex conjugate of both sides of this equation. For the right-hand side, the definition of $f_\r$ in \eqref{eq:fr-def} implies that $\ol{f_{\r'} (-\bar{\tau})} = f_\r (\tau)$ and
	\begin{equation}
		\ol{f_{\{S\r'\}} (S \cdot (-\bar{\tau}))}
		= \ol{f_{\left\{\prmatrix*{-r_2}{-r_1}\right\}} (S \cdot (-\bar{\tau}))}
		= f_{\left\{\prmatrix*{r_2}{-r_1}\right\}} (S \cdot \tau)
		= f_{\{-S\r\}} (S \cdot \tau),
	\end{equation}
	so the right-hand side becomes $f_\r (\tau) - f_{\{-S\r\}} (S \cdot \tau)$. For the left side, we have
	\begin{equation} \label{eq:s0-conj}
		\ol{\S_0 (\tf_{\r'}) (-\bar{\tau})}
		= \int_{0}^{\infty} \ol{e^{\xi/\bar{\tau}} \B(\psi_{\r'}) (\xi)}\, d\xi
		= -\int_{0}^{-\infty} e^{-\xi/\tau} \ol{\B(\psi_{\r'}) (-\xi)}\, d\xi,
	\end{equation}
	where we replaced $\xi$ by $-\xi$ in the last integral. By \zcref{prop:borel-transform-psi-r},
	\begin{align}
		\ol{\B(\psi_{\r'}) (-\xi)}
		&= \frac{1}{2\pi \ol{i}} \sum_{n \ne 0} \frac{1}{n^2} \cdot \frac{\ol{\e(n r_2)}}{1 - \ol{\e(-r_1)} e^{\xi/n}}
		= -\frac{1}{2\pi i} \sum_{n \ne 0} \frac{1}{n^2} \cdot \frac{\e(-n r_2)}{1 - \e(r_1) e^{\xi/n}} \\
		&= -\frac{1}{2\pi i} \sum_{n \ne 0} \frac{1}{n^2} \cdot \frac{\e(n r_2)}{1 - \e(r_1) e^{-\xi/n}}
		= -\B(\psi_\r) (\xi),
	\end{align}
	where we replaced $n$ by $-n$ in the sum. By substituting into \eqref{eq:s0-conj}, we find
	\begin{equation}
		\ol{\S_0 (\tf_{\r'}) (-\bar{\tau})}
		= \int_{0}^{-\infty} e^{-\xi/\tau} \B(\psi_\r) (\xi)\, d\xi
		= \S_\pi (\tf_\r) (\tau).
	\end{equation}
	Finally, taking the complex conjugate of both sides of \eqref{eq:s0-tau-bar} yields
	\begin{equation}
		\S_\pi (\tf_\r) (\tau)
		= f_\r (\tau) - f_{\{-S\r\}} (S \cdot \tau). \tag*{\qedhere}
	\end{equation}
\end{proof}

\begin{remark} \label{rmk:bl-analytic-cont}
	Since the poles of $\B(\psi_\r) (\tau)$ lie on the imaginary axis, a standard contour argument shows that
	\begin{alignat}{3}
		\S_\theta (\tf_\r) (\tau) &= \S_0   (\tf_\r) (\tau) &\quad &\text{for} \quad &-\frac{\pi}{2} &< \theta < \frac{\pi}{2}, \\
		\S_\theta (\tf_\r) (\tau) &= \S_\pi (\tf_\r) (\tau) &\quad &\text{for} \quad & \frac{\pi}{2} &< \theta < \frac{3\pi}{2},
	\end{alignat}
	whenever both sides are defined. But $\B(\psi_\r) (\xi)$ is of exponential type $0$ on the ray $e^{i\theta} \R_+$, so $\S_\theta (\tf_\r) (\tau)$ is analytic on the half-plane $\Re(e^{-i\theta} \tau) > 0$, and hence, by varying $\theta$, the above equations provide an analytic continuation for $\S_0 (\tf_\r) (\tau)$ to the cut plane $\C' = \C_S = \C \sm \R_{\le 0}$ and for $\S_\pi (\tf_\r) (\tau)$ to $\C_{-S} = \C \sm \R_{\ge 0}$.
\end{remark}

\subsubsection{Borel--Laplace sums for \texorpdfstring{$r_1 \in \Z$}{r₁ ∈ ℤ}}
\begin{theorem} \label{thm:bl-fr-r1-int}
	Let $\r = \prmatrix*{r_1}{r_2} \in \Z \times \R_{> 0}$ and $\tau \in \H$. The Borel--Laplace sums of the asymptotic series of $f_\r (\tau)$ as $\tau \to 0$ along the rays $\theta = 0$ and $\theta = \pi$ satisfy
	\begin{alignat}{2}
		\S_0   (\tf_\r) (\tau) &= f_\r (\tau) - f_{\{ S\r\}} (S \cdot \tau), & \qquad &\Re(\tau) > 0, \label{eq:s0-fr-r1-int} \\
		\S_\pi (\tf_\r) (\tau) &= f_\r (\tau) - f_{\{-S\r\}} (S \cdot \tau), & \qquad &\Re(\tau) < 0. \label{eq:s-pi-fr-r1-int}
	\end{alignat}
\end{theorem}

\begin{proof}
	We only prove \eqref{eq:s0-fr-r1-int}, since \eqref{eq:s-pi-fr-r1-int} follows from it by the same argument used to deduce \eqref{eq:s-pi-fr-r1-nonint} from \eqref{eq:s0-fr-r1-nonint} in \zcref{thm:bl-fr-r1-nonint}. Moreover, we will prove \eqref{eq:s0-fr-r1-int} for $r_2 \in (0,1]$, then the result follows for $r_2 > 1$ by induction as in the proof of \zcref{thm:bl-fr-r1-nonint}.
	
	Consider first the case $r_2 = 1$. By noting that $B_k (1) = B_k = 0$ for odd $k \ge 3$, we find, from the definition of $\psi_{r_2}$ stated in \zcref{lem:fr-asymp}, that $\psi_{r_2} (\tau) = -\pi i \tau / 12$. Therefore, by \zcref{lem:fr-asymp,cor:BL-analytic},
	\begin{align}
		\S_0 (\tf_\r) (\tau)
		&= -\frac{\pi i}{12\tau} - B_1 (1) \log(-2\pi i \tau) - \log\frac{\Gamma(1)}{\sqrt{2\pi}} + \S_0 (\psi_{r_2}) (\tau) \\
		&= -\frac{\pi i}{12\tau} - \frac{1}{2} \log(-2\pi i \tau) + \frac{1}{2} \log(2\pi) - \frac{\pi i \tau}{12} \\
		&= -\frac{\pi i}{12\tau} - \frac{1}{2} \log(-i \tau) - \frac{\pi i \tau}{12}.
	\end{align}
	On other hand, since $r_1 \in \Z$, we have $f_\r (\tau) = \slog\qpoch{\e(\tau)}{\e(\tau)}$. Recall that the Dedekind eta function is given by $\eta(\tau) = e^{\pi i\tau/12}\, \qpoch{\e(\tau)}{\e(\tau)}$, so
	\begin{equation}
		\slog(\eta(\tau)) = \frac{\pi i\tau}{12} + f_\r (\tau).
	\end{equation}
	Using its transformation formula
	\begin{equation}
		\eta(-1/\tau) = (-i\tau)^{1/2} \eta(\tau),
	\end{equation}
	we obtain
	\begin{equation}
		-\frac{\pi i}{12\tau} + f_\r (-1/\tau)
		= \frac{1}{2} \log(-i\tau) + \frac{\pi i\tau}{12} + f_\r (\tau).
	\end{equation}
	Rearranging the above equation and writing $-1/\tau = S \cdot \tau$ yields
	\begin{equation}
		f_\r (\tau) - f_\r (S \cdot \tau)
		= -\frac{\pi i}{12\tau} - \frac{1}{2} \log(-i\tau) - \frac{\pi i\tau}{12}
		= S_0 (\tf_\r) (\tau).
	\end{equation}
	But $\{S\r\} = \left\{\prmatrix*{-1}{r_1}\right\} = \prmatrix*{0}{1}$, so periodicity of $f_\r$ in the first parameter implies that $f_\r (S \cdot \tau) = f_{\{S\r\}} (S \cdot \tau)$. This proves \eqref{eq:s0-fr-r1-int} for $r_2 = 1$.
		
	Suppose now that $r_2 \in (0,1)$ and let $0 < \Re(\tau) < 1/r_2$. Set
	\begin{equation}
		K(t) = \frac{1}{1 - e^t} + \frac{1}{t}
		\qqand
		A_\pm (t) = \log(1 - e^{-t/\tau} \e(\pm r_2)).
	\end{equation}
	By \zcref{prop:borel-transform-psi-r2}, we have
	\begin{equation}
		\S_0 (\psi_{r_2}) (\tau)
		= \int_{0}^{\infty} e^{-\xi/\tau} \B(\psi_{r_2}) (\xi)\, d\xi
		= \frac{1}{2\pi i} \int_{0}^{\infty} e^{-\xi/\tau} \sum_{n \ne 0} \frac{\e(n r_2)}{n^2}\, K(-\xi/n)\, d\xi.
	\end{equation}
	Since $K(t)$ is uniformly bounded on $\R$ (the singularity at $t = 0$ is removable), we may interchange the sum and the integral by Fubini's Theorem to obtain
	\begin{equation}
		\S_0 (\psi_{r_2}) (\tau)
		= \frac{1}{2\pi i} \sum_{n=1}^{\infty} \frac{1}{n^2}
		\int_{0}^{\infty} e^{-\xi/\tau} \bigl(\e(n r_2) K(-\xi/n) + \e(-n r_2) K(\xi/n)\bigr)\, d\xi.
	\end{equation}
	By substituting $\xi = nt$ in the integral, the above becomes
	\begin{align}
		\S_0 (\psi_{r_2}) (\tau)
		&= \frac{1}{2\pi i} \sum_{n=1}^{\infty} \frac{1}{n}
		\int_{0}^{\infty} e^{-nt/\tau} \bigl(K(-t) \e(n r_2) + K(t) \e(-n r_2)\bigr)\, dt. \\
		&= \frac{1}{2\pi i} \left(\int_{0}^{\infty} K(-t) \sum_{n=1}^{\infty} \frac{e^{-nt/\tau} \e(n r_2)}{n}\, dt
		+ \int_{0}^{\infty} K(t) \sum_{n=1}^{\infty} \frac{e^{-nt/\tau} \e(-n r_2)}{n}\, dt\right). \\
		&= -\frac{1}{2\pi i} \left(\int_{0}^{\infty} K(-t) A_+ (t)\, dt + \int_{0}^{\infty} K(t) A_- (t)\, dt\right),
	\end{align}
	where used the Taylor expansion of $\log(1 - z)$ in the last step. Substituting the definition of $K(t)$ and combining the integrals gives
	\begin{align}
		\S_0 (\psi_{r_2}) (\tau)
		&= -\frac{1}{2\pi i} \int_{0}^{\infty} \left(\frac{A_+ (t)}{1 - e^{-t}} + \frac{A_- (t)}{1 - e^t} + \frac{A_- (t) - A_+ (t)}{t}\right) dt.
	\end{align}
	
	For $\varepsilon > 0$, let
	\begin{equation}
		I_1 (\varepsilon) = \int_{\varepsilon}^{\infty} \left(\frac{A_+ (t)}{1 - e^{-t}} + \frac{A_- (t)}{1 - e^t}\right) dt
		\qqand
		I_2 (\varepsilon) = \int_{\varepsilon}^{\infty} \frac{A_- (t) - A_+ (t)}{t}\, dt
	\end{equation}
	so that
	\begin{equation} \label{eq:s0-r2-lim}
		\S_0 (\psi_{r_2}) (\tau) = -\frac{1}{2\pi i} \lim_{\varepsilon \to 0^+} (I_1 (\varepsilon) + I_2 (\varepsilon)).
	\end{equation}
	The cutoff at $\varepsilon > 0$ here is necessary to split the integral since both integrals diverge at $\varepsilon = 0$.
	For $t > 0$, we have that
	\begin{equation}
		A_\pm (t)
		= \log(1 - e^{-t/\tau} \e(\pm r_2))
		= \ell(-t/\tau \pm 2\pi i r_2).
	\end{equation}
	Substituting the above into $I_1 (\varepsilon)$ then splitting it into two integrals and replacing $t$ by $-t$ in the first integral yields
	\begin{align}
		I_1 (\varepsilon)
		&= \int_{\varepsilon}^{\infty} \left(\frac{\ell(-t/\tau + 2\pi i r_2)}{1 - e^{-t}} + \frac{\ell(-t/\tau - 2\pi i r_2)}{1 - e^t}\right) dt \\
		&= \int_{-\infty}^{-\varepsilon} \frac{\ell(t/\tau + 2\pi i r_2)}{1 - e^t}\, dt
		+ \int_{\varepsilon}^{\infty} \frac{\ell(-t/\tau - 2\pi i r_2)}{1 - e^t}\, dt.
	\end{align}
	By \eqref{eq:l-reflection}, we have for $t > 0$,
	\begin{equation}
		\ell(-t/\tau - 2\pi i r_2)
		= \ell(t/\tau + 2\pi i r_2) - \frac{t}{\tau} - 2\pi i B_1 (r_2),
	\end{equation}
	so
	\begin{equation}
		I_1 (\varepsilon)
		= \left(\int_{-\infty}^{-\varepsilon} + \int_{\varepsilon}^{\infty}\right) \frac{\ell(t/\tau + 2\pi i r_2)}{1 - e^t}\, dt
		- \int_{\varepsilon}^{\infty} \frac{t/\tau}{1 - e^t}\, dt
		- \int_{\varepsilon}^{\infty} \frac{2\pi i B_1 (r_2)}{1 - e^t}\, dt.
	\end{equation}
	In the limit, the sum of the first two integrals converges to the principal value, and the third integral becomes
	\begin{equation}
		\int_{0}^{\infty} \frac{t/\tau}{1 - e^t}\, dt
		= -\frac{1}{\tau} \Li_2 (1)
		= -\frac{\pi^2}{6\tau}.
	\end{equation}
	The last integral is
	\begin{equation}
		\int_{\varepsilon}^{\infty} \frac{2\pi i B_1 (r_2)}{1 - e^t}\, dt
		= 2\pi i B_1 (r_2) \log(1 - e^{-\varepsilon}),
	\end{equation}
	which diverges as $\varepsilon \to 0^+$. Consequently,
	\begin{equation} \label{eq:lim-I1}
		\lim_{\varepsilon \to 0^+} (I_1 (\varepsilon) + 2\pi i B_1 (r_2) \log(1 - e^{-\varepsilon}))
		= \PV \int_{-\infty}^{\infty} \frac{\ell(t/\tau + 2\pi i r_2)}{1 - e^t}\, dt + \frac{\pi^2}{6\tau}
	\end{equation}
	
	Now consider
	\begin{equation}
		\lim_{\varepsilon \to 0^+} (I_2 (\varepsilon) - 2\pi i B_1 (r_2) \log(1 - e^{-\varepsilon})).
	\end{equation}
	Since $\log(1 - e^{-\varepsilon}) \sim \log\varepsilon$ as $\varepsilon \to 0^+$, this limit equals
	\begin{equation}
		\lim_{\varepsilon \to 0^+} (I_2 (\varepsilon) - 2\pi i B_1 (r_2) \log\varepsilon).
	\end{equation}
	Let
	\begin{equation}
		B(s,t) = \log(1 - e^{-st} \e(-r_2)) - \log(1 - e^{-st} \e(r_2))
	\end{equation}
	so that $A_- (t) - A_+ (t) = B(1/\tau, t)$. Using $\log\varepsilon = \int_{\varepsilon}^{1} 1/t\, dt$, we get
	\begin{align}
		I_2 (\varepsilon) - 2\pi i B_1 (r_2) \log\varepsilon
		&= \int_{\varepsilon}^{\infty} \frac{A_- (t) - A_+ (t)}{t}\, dt
		+ 2\pi i B_1 (r_2) \int_{\varepsilon}^{1} \frac{1}{t}\, dt \\
		&= \int_{\varepsilon}^{1} \frac{B(1/\tau, t) + 2\pi i B_1 (r_2)}{t}\, dt
		+ \int_{1}^{\infty} \frac{B(1/\tau, t)}{t}\, dt,
	\end{align}
	The first integral now has a finite limit as $\varepsilon \to 0^+$ because $A_\pm (t)$ are analytic near $t = 0$ and
	\begin{equation} \label{eq:B0}
		B(1/\tau, 0)
		= \log(1 - \e(-r_2)) - \log(1 - \e(r_2))
		= 2i \Im\log(1 - \e(-r_2))
		= -2\pi i B_1 (r_2),
	\end{equation}
	where we write $1 - \e(-r_2) = 2\sin(\pi r_2) e^{-\pi i B_1 (r_2)}$ to get the last equality. It follows that
	\begin{equation} \label{eq:B-asymptotic}
		B(1/\tau, t) + 2\pi i B_1 (r_2) = O(t) \quad \text{as } t \to 0,
	\end{equation}
	so
	\begin{equation} \label{eq:lim-I2}
		\lim_{\varepsilon \to 0^+} (I_2 (\varepsilon) - 2\pi i B_1 (r_2) \log\varepsilon)
		= \int_{0}^{1} \frac{B(1/\tau, t) + 2\pi i B_1 (r_2)}{t}\, dt + \int_{1}^{\infty} \frac{B(1/\tau, t)}{t}\, dt.
	\end{equation}
	
	For $\Re(s) > 0$, define
	\begin{equation}
		F(s) = \int_{0}^{1} \frac{B(s, t) + 2\pi i B_1 (r_2)}{t}\, dt + \int_{1}^{\infty} \frac{B(s, t)}{t}\, dt.
	\end{equation}
	Our goal is to get a closed-form expression of $F(s)$ then evaluate it at $s = 1/\tau$. To find such an expression, we differentiate under the integral sign. This is possible by the Dominated Convergence Theorem because both integrands have the derivative
	\begin{equation}
		\frac{1}{t} \cdot \frac{\pd}{\pd s}\, B(s,t)
		= \frac{1}{t} \left(\frac{t e^{-st} \e(-r_2)}{1 - e^{-st} \e(-r_2)} - \frac{t e^{-st} \e(r_2)}{1 - e^{-st} \e(r_2)}\right)
		= \frac{e^{-st} \e(-r_2)}{1 - e^{-st} \e(-r_2)} - \frac{e^{-st} \e(r_2)}{1 - e^{-st} \e(r_2)},
	\end{equation}
	which is uniformly bounded for $t \in [0,1]$ by
	\begin{equation}
		\frac{\abs{1 + \e(r_2)}}{\abs{1 - \e(r_2)}},
	\end{equation}
	and for $t \ge 1$ and $\Re(s) \ge \delta > 0$,
	\begin{equation}
		\abs*{\frac{1}{t} \cdot \frac{\pd}{\pd s}\, B(s,t)}
		\le \frac{\abs{e^{-st}}}{\abs{1 - e^{-st} \e(-r_2)}} + \frac{\abs{e^{-st}}}{\abs{1 - e^{-st} \e(r_2)}}
		\le \frac{2\abs{e^{-st}}}{1 - \abs{e^{-st}}}
		= \frac{2 e^{-\delta t}}{1 - e^{-\delta t}},
	\end{equation}
	which is integrable on $[1,\infty)$, Therefore,
	\begin{align}
		F'(s)
		&= \int_{0}^{1} \frac{1}{t}\cdot \frac{\pd}{\pd s} (B(s, t) + 2\pi i B_1 (r_2))\, dt
		+ \int_{1}^{\infty} \frac{1}{t} \cdot \frac{\pd}{\pd s}\, B(s, t)\, dt \\
		&= \int_{0}^{\infty} \frac{1}{t} \cdot \frac{\pd}{\pd s}\, B(s, t)\, dt \\
		&= \frac{1}{s} \left(\log(1 - e^{-st} \e(-r_2)) - \log(1 - e^{-st} \e(r_2))\right)\big|_{t=0}^{\infty} \\
		&= -\frac{1}{s} (\log(1 - \e(-r_2)) - \log(1 - \e(r_2))) \\
		&= \frac{2\pi i B_1 (r_2)}{s},
	\end{align}
	where we used \eqref{eq:B0} for the last equality. Integrating with respect to $s$ yields
	\begin{equation} \label{eq:F-C}
		F(s) = 2\pi i B_1 (r_2) \log s + C
	\end{equation}
	for some constant $C$. Setting $s = 1$ yields $C = F(1)$.
	
	To evaluate $F(1)$, we use the series expansion
	\begin{align}
		B(1,t)
		&= \log(1 - e^{-t} \e(-r_2)) - \log(1 - e^{-t} \e(r_2)) \\
		&= -\sum_{n=1}^{\infty} \frac{e^{-nt} \e(-n r_2)}{n} + \sum_{n=1}^{\infty} \frac{e^{-nt} \e(n r_2)}{n} \\
		&= 2i \sum_{n=1}^{\infty} \frac{e^{-nt} \sin(2\pi n r_2)}{n}
	\end{align}
	and the Fourier expansion of $B_1 (x)$
	\begin{equation} \label{eq:B1-fourier}
		B_1 (r_2)
		= -\frac{1}{\pi} \sum_{n=1}^{\infty} \frac{\sin(2\pi n r_2)}{n}
	\end{equation}
	to get
	\begin{align}
		F(1)
		&= \int_{0}^{1} \frac{B(1, t) + 2\pi i B_1 (r_2)}{t}\, dt + \int_{1}^{\infty} \frac{B(1, t)}{t}\, dt \\
		&= \int_{0}^{1} \frac{2i}{t} \left(\sum_{n=1}^{\infty} \frac{e^{-nt} \sin(2\pi n r_2)}{n} - \sum_{n=1}^{\infty} \frac{\sin(2\pi n r_2)}{n}\right)\, dt
		+ \int_{1}^{\infty} \frac{2i}{t} \sum_{n=1}^{\infty} \frac{e^{-nt} \sin(2\pi n r_2)}{n}\, dt \\
		&= 2i \sum_{n=1}^{\infty} \frac{\sin(2\pi n r_2)}{n} \int_{0}^{1} \frac{e^{-nt} - 1}{t}\, dt
		+ 2i \sum_{n=1}^{\infty} \frac{\sin(2\pi n r_2)}{n} \int_{1}^{\infty} \frac{e^{-nt}}{t}\, dt \\
		&= 2i \sum_{n=1}^{\infty} \frac{\sin(2\pi n r_2)}{n} \left(\int_{0}^{1} \frac{e^{-nt} - 1}{t}\, dt
		+ \int_{1}^{\infty} \frac{e^{-nt}}{t}\, dt\right). \label{eq:F1-series}
	\end{align}
	The second integral in the above is the exponential integral $E_1 (n)$, which has the series expansion
	\begin{equation}
		E_1 (n) = -\gamma - \log n - \sum_{k=1}^{\infty} \frac{(-n)^k}{k \cdot k!},
	\end{equation}
	and the first integral can be evaluated by substituting the series expansion of $e^{-nt}$ as
	\begin{equation}
		\int_{0}^{1} \frac{e^{-nt} - 1}{t}\, dt
		= \int_{0}^{1} \sum_{k=1}^{\infty} \frac{(-n)^k t^{k - 1}}{k!}\, dt
		= \sum_{k=1}^{\infty} \frac{(-n)^k}{k \cdot k!}
	\end{equation}
	Substituting the above two equations into \eqref{eq:F1-series} then using \eqref{eq:B1-fourier} again yields
	\begin{equation}
		F(1)
		= 2i \sum_{n=1}^{\infty} \frac{\sin(2\pi n r_2)}{n} (-\gamma - \log n)
		= 2\pi i \gamma B_1 (r_2) - 2i \sum_{n=1}^{\infty} \frac{\sin(2\pi n r_2)}{n} \log n.
	\end{equation}
	Finally, by substituting $x = r_2$ into Kummer's Fourier expansion of $\log\Gamma(x)$ \cite[Thm.~1.7.1]{andrews1999special},
	\begin{equation*}
		\log\frac{\Gamma(x)}{\sqrt{2\pi}}
		= -\frac{1}{2} \log(2\sin(\pi x)) - (\gamma + \log(2\pi)) \left(x - \frac{1}{2}\right)
		+ \frac{1}{\pi} \sum_{n=1}^{\infty} \frac{\sin(2\pi n x)}{n} \log n, \quad 0 < x < 1,
	\end{equation*}
	where $\gamma$ is the Euler--Mascheroni constant, we deduce that
	\begin{equation}
		F(1)
		= -\pi i \log(2\sin(\pi r_2)) - 2\pi i B_1 (r_2) \log(2\pi) - 2\pi i \log\frac{\Gamma(r_2)}{\sqrt{2\pi}}.
	\end{equation}
	Taking log of both sides of $1 - \e(r_2) = 2\sin(\pi r_2) e^{\pi i B_1 (r_2)}$ and rearranging yields
	\begin{equation}
		\log(2\sin(\pi r_2)) = \log(1 - \e(r_2)) - \pi i B_1 (r_2),
	\end{equation}
	Moreover,
	\begin{equation}
		\log(2\pi)
		= \log(i (-2\pi i))
		= \log(i) + \log(-2\pi i)
		= \frac{\pi i}{2} + \log(-2\pi i).
	\end{equation}
	Therefore,
	\begin{equation}
		F(1)
		= -\pi i \log(1 - \e(r_2)) - 2\pi i B_1 (r_2) \log(-2\pi i) - 2\pi i \log\frac{\Gamma(r_2)}{\sqrt{2\pi}}.
	\end{equation}
	
	Substituting the above value for $C = F(1)$ into \eqref{eq:F-C} yields
	\begin{equation}
		F(s) = 2\pi i B_1 (r_2) \log s - \pi i \log(1 - \e(r_2)) - 2\pi i B_1 (r_2) \log(-2\pi i) - 2\pi i \log\frac{\Gamma(r_2)}{\sqrt{2\pi}}.
	\end{equation}
	Recall that the limit which was to be evaluated \eqref{eq:lim-I2} is precisely $F(1/\tau)$, which equals
	\begin{align}
		F(1/\tau)
		&= 2\pi i B_1 (r_2) \log(1/\tau) - \pi i \log(1 - \e(r_2)) - 2\pi i B_1 (r_2) \log(-2\pi i)
		- 2\pi i \log\frac{\Gamma(r_2)}{\sqrt{2\pi}} \\
		&= -\pi i \log(1 - \e(r_2)) - 2\pi i B_1 (r_2) \log(-2\pi i \tau) - 2\pi i \log\frac{\Gamma(r_2)}{\sqrt{2\pi}},
	\end{align}
	where we used $0 < \arg(\tau) < \pi/2$ to manipulate the logarithms. Therefore, by summing \eqref{eq:lim-I1} and \eqref{eq:lim-I2} then substituting $F(1/\tau)$ for the second limit, we find
	\begin{align}
		\lim_{\varepsilon \to 0^+} (I_1 (\varepsilon) + I_2 (\varepsilon))
		&= \lim_{\varepsilon \to 0^+} (I_1 (\varepsilon) + 2\pi i B_1 (r_2) \log(1 - e^{-\varepsilon}))
		+ \lim_{\varepsilon \to 0^+} (I_2 (\varepsilon) - 2\pi i B_1 (r_2) \log\varepsilon) \\
		&=
		\begin{multlined}[t]
			\PV \int_{-\infty}^{\infty} \frac{\ell(t/\tau + 2\pi i r_2)}{1 - e^t}\, dt + \frac{\pi^2}{6\tau} - \pi i \log(1 - \e(r_2)) \\
			- 2\pi i B_1 (r_2) \log(-2\pi i \tau) - 2\pi i \log\frac{\Gamma(r_2)}{\sqrt{2\pi}}
		\end{multlined}
	\end{align}
	
	Recall that
	\begin{equation}
		\tf_\r (\tau)
		= -\frac{\pi i}{12\tau} - B_1 (r_2) \log(-2\pi i \tau) - \log\frac{\Gamma(r_2)}{\sqrt{2\pi}} + \psi_{r_2} (\tau)
	\end{equation}
	and
	\begin{equation}
		\S_0 (\psi_{r_2}) (\tau) = -\frac{1}{2\pi i} \lim_{\varepsilon \to 0^+} (I_1 (\varepsilon) + I_2 (\varepsilon)).
	\end{equation}
	Applying \zcref{cor:BL-analytic} then using \zcref{rmk:int-cont-shift} therefore yields
	\begin{align}
		\S_0 (\tf_\r) (\tau)
		&= -\frac{\pi i}{12\tau} - B_1 (r_2) \log(-2\pi i \tau) - \log\frac{\Gamma(r_2)}{\sqrt{2\pi}} + \S_0 (\psi_{r_2}) (\tau) \\
		&= -\frac{1}{2\pi i} \PV \int_{-\infty}^{\infty} \frac{\ell(t/\tau + 2\pi i r_2)}{1 - e^t}\, dt + \frac{1}{2} \log(1 - \e(r_2)) \\
		&= -\frac{1}{2\pi i} \int_{-\infty}^{\infty} \frac{\ell(t/\tau)}{1 - \e(-r_2 \tau) e^t}\, dt.
	\end{align}
	Finally, recalling that $0 < \Re(\tau) < 1/r_2$ and using \zcref{cor:int-to-qpoch} with $z = -r_2 \tau$ yields
	\begin{equation}
		\S_0 (\tf_\r) (\tau) = f_\r (\tau) - f_{\{S\r\}} (S \cdot \tau).
	\end{equation}
	Since both sides of the above equation are analytic for $\Re(\tau) > 0$, this proves \eqref{eq:s0-fr-r1-int} for $r_2 \in (0,1)$.
\end{proof}

\begin{remark} \label{rmk:FR-conj}
	We observe that Conjecture~5 of Fantini--Rella \cite{fantini2025modular} is equivalent to \zcref{thm:bl-fr-r1-int}. Indeed, if $\r = \prmatrix*{0}{k/N}$ with $N \ge 2$ and $1 \le k \le N - 1$, then
	\begin{align}
		g_{k,N} (\tau)              &= f_\r (N\tau),                                                                       \\
		f_{k,N} (\tau)              &= f_{ S\r} (\tau) = f_{\{ S\r\}} (\tau) - \log(1 - e^{ 2\pi i k/N}), \quad \text{and} \\
		f_{\underline{k} ,N} (\tau) &= f_{-S\r} (\tau) = f_{\{-S\r\}} (\tau) - \log(1 - e^{-2\pi i k/N}),
	\end{align}
	where $g_{k,N}$ and $f_{k,N}$ are as defined in \cite[Eq.~(1.3)]{fantini2025modular} and $\underline{k} \coloneqq N - k$. It now follows from \zcref{thm:bl-fr-r1-int} that
	\begin{alignat}{2}
		\S_0 (\tilde{g}_{k,N}) (\tau)
		&= g_{k,N} (\tau) - f_{k,N} \left(-\frac{1}{N\tau}\right) + \log(1 - e^{2\pi i k/N}), \qquad &\Re(\tau) > 0, \\
		\S_\pi (\tilde{g}_{k,N}) (\tau)
		&= g_{k,N} (\tau) - f_{\underline{k}, N} \left(-\frac{1}{N\tau}\right) + \log(1 - e^{-2\pi i k/N}), \qquad &\Re(\tau) < 0.
	\end{alignat}
	These are precisely the two identities in \cite[Conj.~5]{fantini2025modular}.
\end{remark}

\subsection{Median resummation} \label{sec:med-resum}
By \zcref{rmk:bl-analytic-cont}, the Borel--Laplace sums in \zcref{thm:bl-fr} admit analytic continuations from their respective domains to $\H$. Averaging these continuations as in \eqref{eq:med-resum-disc} gives the following formula for the median resummation.

\begin{corollary} \label{cor:fr-med-resum}
	\frmedresum
\end{corollary}

Notice that the term
\begin{equation}
	f_{\{S\r\}} (S \cdot \tau) + f_{\{-S\r\}} (S \cdot \tau)
\end{equation}
appearing in the median resummation of $\tf_\r$ is symmetric under the involution $\r \mapsto -\r$. This motivates the following result.

\begin{corollary} \label{cor:Fr-med-resum}
	\Frmedresum
\end{corollary}

\begin{proof}
	Throughout the proof, we use the identity
	\begin{equation}
		\{S\{-\r\}\} = \{-S\r\} = \{-\{S\r\}\},
	\end{equation}
	which follow from $S\{-\r\} \equiv -S\r \equiv -\{S\r\} \pmod{\Z^2}$.
	
	Let $\Re(\tau) > 0$. By linearity of the Borel--Laplace sum and \zcref{thm:bl-fr}, we have
	\begin{equation} \label{eq:s0-Fr}
		\begin{aligned}
			\S_0 (\tF_\r) (\tau)
			&= \S_0 (\tf_\r) (\tau) - \S_0 (\tf_{\{-\r\}}) (\tau) \\
			&= \bigl(f_\r (\tau) - f_{\{S\r\}} (S \cdot \tau)\bigr) - \bigl(f_{\{-\r\}} (\tau) - f_{\{S\{-\r\}\}} (S \cdot \tau)\bigr) \\
			&= \bigl(f_\r (\tau) - f_{\{-\r\}} (\tau)\bigr) - \bigl(f_{\{S\r\}} (S \cdot \tau) - f_{\{-S\r\}} (S \cdot \tau)\bigr) \\
			&= F_\r (\tau) - F_{\{S\r\}} (S \cdot \tau).
		\end{aligned}
	\end{equation}
	Similarly, \zcref{thm:bl-fr} implies that for $\Re(\tau) < 0$,
	\begin{equation}
		\begin{aligned}
			\S_\pi (\tF_\r) (\tau)
			&= \S_\pi (\tf_\r) (\tau) - \S_\pi (\tf_{\{-\r\}}) (\tau) \\
			&= \bigl(f_\r (\tau) - f_{\{-S\r\}} (S \cdot \tau)\bigr) - \bigl(f_{\{-\r\}} (\tau) - f_{\{-S\{-\r\}\}} (S \cdot \tau)\bigr) \\
			&= \bigl(f_\r (\tau) - f_{\{-\r\}} (\tau)\bigr) - \bigl(f_{\{-S\r\}} (S \cdot \tau) - f_{\{S\r\}} (S \cdot \tau)\bigr) \\
			&= F_\r (\tau) - F_{\{-S\r\}} (S \cdot \tau).
		\end{aligned}
	\end{equation}
	But $F_{\{-S\r\}} (S \cdot \tau) = F_{\{-S\r\}} (S \cdot \tau)$, so we may also write
	\begin{equation} \label{eq:s-pi-Fr}
		\S_\pi (\tF_\r) (\tau)
		= F_\r (\tau) + F_{\{S\r\}} (S \cdot \tau).
	\end{equation}
	By deforming the contour of integration, we can analytically continue $\S_0 (\tF_\r)$ and $\S_\pi (\tF_\r)$ to the upper-half plane (see \zcref{rmk:bl-analytic-cont}). Finally, averaging \eqref{eq:s0-Fr} and \eqref{eq:s-pi-Fr} yields
	\begin{equation}
		\Smed_{\pi/2} (\tF_\r) (\tau) = F_\r (\tau). \tag*{\qedhere}
	\end{equation}
\end{proof}

\section{Quantum modularity of \texorpdfstring{$f_\r$}{fᵣ}} \label{sec:fr-quantum}
In this section, we will use \zcref{thm:bl-fr,rmk:bl-analytic-cont} to establish that if $\r \in (\R \times [0,1]) \sm (\Z \times \{0\})$, then $f_\r$ is a holomorphic quantum modular form for the group
\begin{equation}
	\Gamma_\r \coloneqq \{\gamma \in \SL_2 (\Z) \colon \gamma \r - \r \in \Z^2\},
\end{equation}
which is a congruence subgroup of $\SL_2 (\Z)$ if $\r \in \Q^2$. In particular, if $N$ is the least common denominator of $r_1$ and $r_2$, then $\Gamma(N) \subseteq \Gamma_\r$ where $\Gamma(N)$ is the principal congruence subgroup of level $N$. For $\gamma \in \Gamma_\r$, define
\begin{equation} \label{eq:fr-cocycle}
	h_\gamma (\tau) = f_\r (\tau) - f_\r (\gamma \cdot \tau).
\end{equation}
Our goal is to show that $h_\gamma$ extends to the cut plane $\C_\gamma$ (defined in \zcref{def:holomorphic-qmf}), hence proving that $f_\r$ is a holomorphic quantum modular form for $\Gamma_\r$.

To this end, we consider
\begin{equation}
	H_\gamma (\r; \tau) \coloneqq f_\r (\tau) - f_{\{\gamma\r\}} (\gamma \cdot \tau).
\end{equation}
We will show that $H_\gamma$ continues analytically to $\C_\gamma$ for $\gamma \in \SL_2 (\Z)$. In fact, this suggests viewing $f_\r (\tau)$ as a function of both $\r$ and $\tau$ and extending the slash action to act on both variables. We therefore write $f(\r; \tau) \coloneqq f_\r (\tau)$ and define
\begin{equation}
	(f | \gamma) (\r; \tau) \coloneqq f(\{\gamma\r\}; \gamma \cdot \tau)
\end{equation}
so that $H_\gamma = f - f | \gamma$. Then we show, according to this definition of the slash action, that $f$ that is a holomorphic quantum modular form for the full modular group $\SL_2 (\Z)$.

Observe that $H$ satisfies the \emph{cocycle condition}: for $\gamma_1, \gamma_2 \in \SL_2 (\Z)$,
\begin{align}
	H_{\gamma_1 \gamma_2} (\r; \tau)
	&= f_\r (\tau) - f(\{\gamma_1 \gamma_2\r\}; \gamma_1 \gamma_2 \cdot \tau) \\
	&= f_\r (\tau) - f(\{\gamma_2\r\}; \gamma_2 \cdot \tau) \label{eq:H-cocycle}
	+ f(\{\gamma_2\r\}; \gamma_2 \cdot \tau) - f(\{\gamma_1 \gamma_2\r\}; \gamma_1 \gamma_2 \cdot \tau) \\
	&= H_{\gamma_2} (\r; \tau) + H_{\gamma_1} (\{\gamma_2\r\}; \gamma_2 \cdot \tau).
\end{align}
Moreover, by \zcref{thm:bl-fr}, we have
\begin{equation} \label{eq:H-BL}
	H_S (\r; \tau) = \S_0 (\tf_\r) (\tau)
	\qand
	H_{-S} (\r; \tau) = \S_\pi (\tf_\r) (\tau).
\end{equation}

\begin{theorem} \label{thm:f-qm}
	Let $\r = \prmatrix*{r_1}{r_2} \in (\R \times [0,1]) \sm (\Z \times \{0\})$. Then $f$ is a weight-$0$ holomorphic quantum modular form for $\SL_2 (\Z)$.
\end{theorem}

\begin{proof}
	Let $\gamma = \prmatrix*{a & b}{c & d} \in \SL_2 (\Z)$. We proceed by induction on $c$. If $c = 0$, then $\gamma = \pm T^n$ for some $n \in \Z$ (in particular, $n = \pm b$). For $\gamma = T^n$, we have
	\begin{equation}
		\{T^n \r\}
		= \left\{\prmatrix{r_1 + n r_2}{r_2}\right\}
		= \prmatrix{r_1 + n r_2 - \floor{r_1 + n r_2}}{r_2 - \ceil{r_2} + 1}.
	\end{equation}
	Since
	\begin{equation}
		\e((r_2 - \ceil{r_2} + 1) (\tau + n) - r_1 - n r_2 + \floor{r_1 + n r_2})
		= \e((r_2 - \ceil{r_2} + 1) \tau - r_1),
	\end{equation}
	it follows by definition of $f$ that
	\begin{equation}
		f(\{T^n \r\}; T^n \cdot \tau)
		= f(\{\r\}; \tau)
		=
		\begin{cases}
			f(\r; \tau),                      & r_2 \in (0,1], \\
			f(\r; \tau) - \log(1 - \e(-r_1)), & r_2 = 0.
		\end{cases}
	\end{equation}
	Hence,
	\begin{equation}
		H_{T^n} (\r; \tau)
		= f(\r; \tau) - f(\{\r\}; \tau)
		=
		\begin{cases}
			0,                  & r_2 \in (0,1], \\
			\log(1 - \e(-r_1)), & r_2 = 0,
		\end{cases}
	\end{equation}
	which clearly continues to $\C_\gamma = \C_{T^n} = \C$. On the other hand, if $\gamma = -T^n = S^2 T^n$, then applying \eqref{eq:H-cocycle} twice yields
	\begin{align}
		H_{S^2 T^n} (\r; \tau)
		&= H_S (\{ST^n \r\}; ST^n \cdot \tau) + H_S (\{T^n\r\}; T^n \cdot \tau) + H_{T^n} (\r; \tau) \\
		&= H_S (\{ST^n \r\}; -1/(\tau + n)) + H_S (\{T^n\r\}; \tau + n) + H_{T^n} (\r; \tau).
	\end{align}
	By \eqref{eq:H-BL} and \zcref{rmk:bl-analytic-cont}, the first two terms are analytic on $\C \sm [-n,\infty)$ and $\C \sm (-\infty, -n]$ respectively, and the last term is entire by the case $\gamma = T^n$. Therefore, $H_\gamma = H_{S^2 T^n}$ is analytic on $\C \sm \R = \C_{-T^n}$.
	
	Suppose now that $c > 0$ and assume the result holds for any $\gamma' \in \SL_2 (\Z)$ with $0 \le c_{\gamma'} < c$. Write
	\begin{equation}
		a = q c - c',
		\qquad
		0 \le c' < c.
	\end{equation}
	Then $\gamma = T^q S \gamma'$, where
	\begin{equation}
		\gamma' = \prmatrix{c & d}{c' & d'}, \qquad d' = qd - b.
	\end{equation}
	By applying \eqref{eq:H-cocycle} twice, we get
	\begin{equation}
		H_\gamma (\r; \tau)
		=  H_{T^q} (\{S\gamma'\r\}; S\gamma' \cdot \tau) + H_S (\{\gamma'\r\}; \gamma' \cdot \tau) + H_{\gamma'} (\r; \tau).
	\end{equation}
	We claim that all terms on the right-hand side are analytic on $\C_\gamma$. Indeed, the first term vanishes by the $c = 0$ case, and the second term is
	\begin{equation}
		H_S (\{\gamma'\r\}; \gamma' \cdot \tau)
		= \S_0 (\tf_{\{\gamma'\r\}}) (\gamma' \cdot \tau),
	\end{equation}
	which is analytic for $\gamma' \cdot \tau \in \C' = \C \sm \R_{\le 0}$ by \eqref{eq:H-BL} and \zcref{rmk:bl-analytic-cont}. To prove that it is analytic for $\tau \in \C_\gamma$, it suffices show that $\gamma' \cdot \C_\gamma \subseteq \C'$. Let $\tau \in \C_\gamma$. If $\tau$ is non-real, then $\gamma' \cdot \tau$ is also non-real and hence $\gamma' \cdot \tau \in \C'$. If $\tau$ is real, then $\tau \in (-d/c, \infty)$ and
	\begin{equation}
		\gamma' \cdot \tau
		= \frac{c\tau + d}{c'\tau + d'}.
	\end{equation}
	Since $c\tau + d > 0$, we must show that $c'\tau + d' > 0$. Observe that
	\begin{equation} \label{eq:induction-inequality}
		\frac{d'}{c'} - \frac{d}{c}
		= \frac{1}{cc'} \det(\gamma')
		= \frac{1}{cc'}
		> 0,
	\end{equation}
	so
	\begin{equation}
		\tau > -\frac{d}{c} > -\frac{d'}{c'}
	\end{equation}
	and $c'\tau + d' > 0$. Hence $\gamma' \cdot \tau \in \C'$. This shows that $H_S (\{\gamma'\r\}; \gamma' \cdot \tau)$ is analytic on $\C_\gamma$.
	
	Now consider $H_{\gamma'} (\r; \tau)$, which is analytic on $\C_{\gamma'}$ by the induction hypothesis. By \eqref{eq:induction-inequality}, we have $\C_\gamma \subseteq \C_{\gamma'}$. This proves that $H_\gamma$ is analytic on $\C_\gamma$.
	
	It remains to consider the case $c < 0$. Since
	\begin{equation}
		0
		= H_I (\r; \tau)
		= H_{\gamma^{-1} \gamma} (\r; \tau)
		= H_{\gamma^{-1}} (\{\gamma\r\}; \gamma \cdot \tau) + H_\gamma (\r; \tau),
	\end{equation}
	it follows that
	\begin{equation}
		H_\gamma (\r; \tau)
		= -H_{\gamma^{-1}} (\{\gamma\r\}; \gamma \cdot \tau).
	\end{equation}
	But
	\begin{equation}
		\gamma^{-1} = \prmatrix{d & -b}{-c & a},
	\end{equation}
	so the case $c > 0$ implies that $H_{\gamma^{-1}} (\{\gamma\r\}; \gamma \cdot \tau)$ is analytic for $\gamma \cdot \tau \in \C_{\gamma^{-1}}$. Since
	\begin{equation}
		\gamma \cdot \tau \in \C_{\gamma^{-1}}
		\iff -c (\gamma \cdot \tau) + a = \frac{1}{c\tau + d} \in \C'
		\iff c\tau + d \in \C'
		\iff \tau \in \C_\gamma,
	\end{equation}
	we conclude that $H_{\gamma^{-1}} (\{\gamma\r\}; \gamma \cdot \tau) = -H_\gamma (\r; \tau)$ is analytic for $\tau \in \C_\gamma$. This completes the proof.
\end{proof}

\begin{corollary} \label{cor:fr-qm}
	Let $\r = \prmatrix*{r_1}{r_2} \in \R \times [0,1] \sm \Z \times \{0\}$. Then $f_\r$ is a weight-$0$ holomorphic quantum modular form for $\Gamma_\r$.
\end{corollary}

\begin{proof}
	Since $r_2 \in [0,1]$, it follows that
	\begin{equation}
		f_\r (\tau) - f_{\{\r\}} (\tau)
		=
		\begin{cases}
			0,                  & r_2 \in (0,1], \\
			\log(1 - \e(-r_1)), & r_2 = 0.
		\end{cases}
	\end{equation}
	Let $\gamma \in \Gamma_\r$. Then $\gamma \r \equiv \r \pmod{\Z^2}$ and
	\begin{equation}
		H_\gamma (\r; \tau) - h_\gamma (\tau)
		= f_\r (\tau) - f_{\{\gamma\r\}} (\gamma \cdot \tau) - f_\r (\tau) + f_\r (\gamma \cdot \tau)
		= f_\r (\gamma \cdot \tau) - f_{\{\r\}} (\gamma \cdot \tau)
	\end{equation}
	is a constant for fixed $\r$. Since $H_\gamma$ is analytic on $\C_\gamma$ (by \zcref{thm:f-qm}), it follows that $h_\gamma$ is also analytic on $\C_\gamma$, which proves the claim.
\end{proof}

\begin{remark}
	We note that, for $\r \in \Q^2$, Kopp's \cite[Thm.~4.30]{kopp2024shintani} establishes the multiplicative modularity of $\varpi_\r (\tau) = \exp(f_\r (\tau))$ with respect to the ``shin'' cocycle $\shin_\gamma^\r$. In particular, the cocycle associated with $f_\r$ is related to $\shin^\r$ by
	\begin{equation}
		h_\gamma (\tau) = -\slog\shin_\gamma^\r (\tau).
	\end{equation}
\end{remark}

\bibliographystyle{abbrvurl}
\bibliography{references}
\end{document}